# m-Order Time Optimal Control Synthesis Function of Discrete System and Its Application

Qianghui Xiao[a]

[a] *College of Electrical and Information Engineering, Hunan University of Technology, Zhuzhou, Hunan, 412007, PR China*
*E-mail address: 10010@hut.edu.cn;qh.xiao004@163.com*

**Abstract**

In this paper, first of all, we introduce the basic concepts of generating function in combinatorics and some combinatorial identities. In order to facilitate the understanding of $m(m \geq 2)$ -order time optimal control synthesis function of discrete system (referred as m-order synthesis function), secondly, we introduce the derivation process and control ideas of 2nd-order synthesis function $fsun()$, and then deduce in detail the m-order synthesis function by means of generating function. By use of the m-order tracking-form synthesis function with filter factor, the methods of signal extraction and its predictive compensation are presented in this paper, and their immunity and effectiveness are verified by numerical simulation. In the back part of this paper, we analyze the features of various components of 2nd-order active disturbance rejection control (ADRC) technology, and then by constructing the new type structures of extended state observer (ESO) and feedback controller which satisfy the separation principle, we construct a new active disturbance rejection control theory (ADRCT) based on the m-order tracking-form synthesis function with filter factor. Finally, the numerical simulation shows that the new ADRCT has good immunity and robustness.

## 1. Introduction

Qian Xuesen and Song Jian studied the time optimal control problem of discrete system in (Qian & Song, 2011), however, no time optimal control synthesis function of discrete system (referred as synthesis function) based on the state variables (or feedback signals) is given in (Qian & Song, 2011). In the process of proposing and developing the 2nd-order active disturbance rejection control (ADRC) technology, Han Jingqing deduced the approximate 2nd-order synthesis function $fhan()$ by means of the parabola (Han, 2008; Han & Yuan, 1999), and pioneered the introduction of filtering factor to suppress the high-frequency vibration in time optimal control that has been successfully applied to the discrete-time tracking differentiator (Han & Yuan, 1999), signal digital filtering (Han & Huang, 2003; Song, Gan, & Han, 2003; Wu, Lin, & Han, 2004), transition process arrangement (Han, 1998), error feedback control (Han, 1995b), high-frequency vibration elimination (Huang & Han, 1999) and so on. It can be said that the function $fhan()$ is an integral part of ADRC technology. And Sun Biao and Sun Xiuxia found that the function $fhan()$ is not the 2nd-order synthesis function but the likeness of 2nd-order synthesis function of discrete system, the real 2nd-order synthesis function of discrete system, named function $fsun()$, is deduced in (Sun & Sun, 2010). Up to now, there is no research on the $m(m \geq 3)$ -order synthesis function of discrete system. But the latter studies in this paper show that in order to construct the new type structure of extended state observer (ESO) satisfying the separation principle for 2nd-order discrete system, the $m(m \geq 3)$ -order synthesis function needs to be used; meanwhile, the $m(m \geq 3)$ -order and higher order synthesis functions are also required for the $m(m \geq 3)$ -order discrete system. And the above results provide a good starting point for the study of $m(m \geq 3)$ -order synthesis function of discrete system.

As an application of 2nd-order synthesis function, ADRC technology has experienced more than 30 years of tortuous process from the beginning of its design ideas. It has also been 20 years since ADRC was formally put forward, and from its development to today, ADRC technology is known around the world for its unique control ideas and superior control quality (Alonge, et al., 2017a, 2017b; Feng & Guo, 2017; Ran, Wang, & Dong, 2016, 2017; Yao

& Deng, 2017). In the process of its development, through the fully understanding and analyzing the advantages and disadvantages of classical PID and modern control theory, absorbing the essence of classical PID, and combining the modern control theory's understanding of control system and modern signal processing technology, a series of new practical controllers are gradually constructed (Gao, 2013; Han, 2008; Huang & Zhang, 2002).

ADRC technology is formally formed into two stages: 1) the generation of tracking differentiator (Han & Huang, 2003; Han & Yuan, 1999; Song, Gan, & Han, 2003; Wu, Lin, & Han, 2004) and the emergence of a new type of nonlinear PID (Han, 1994a, 1994b); 2) the invention of extended state observer (ESO) (Han, 1995a), which marks the birth of active disturbance rejection controller. Based on the above basic results, a new type of nonlinear error feedback control law that can automatically compensate the system's internal and external disturbances is proposed in (Han, 1995b) and therefore is named as active disturbance rejection controller (Han, 1998).

Since the production and verification of above ideas mostly come from a large number of computer simulation and calculation techniques, rather than the traditional theoretical derivation, people are accustomed to call it ADRC technology instead of active disturbance rejection control theory (ADRCT) (Huang & Zhang, 2002). By investigating its reason, the author of this paper thinks that there are mainly the following questions:

1) So far it is not rigorously proved in theory that the parameter designs of extended state observer and feedback controller satisfy the separation principle, which is one of main reasons why people do not consider it as a theory.

2) Is it appropriate to use the error feedback control law? For the linear error feedback control system, this results in a relative order of one for closed-loop control system, regardless of the relative order of open-loop system.

To analyze the above questions and deal with them appropriately, we need to apply the m-order tracking-form synthesis function with filter factor derived in the front part of this paper. A theoretical analysis in the linear case is given and the appropriate theoretical solutions are put forward in the latter part of this paper.

This paper first introduces the basic concepts of ordinary generating function (referred as generating function) in combinatorics to be used latter, then derives the 2nd-order synthesis function $fsun()$, deduces in detail the $m(m \geq 2)$-order synthesis function $fxiao()$ of discrete system by means of generating function, introduces the filter factor to the $m(m \geq 2)$-order synthesis function, and gives the m-order synthesis function based on the tracking problem. By use of the m-order tracking-form synthesis function with filter factor, it is convenient to extract the filtered feedback signal and its differential signals from 1st-order to $(m-1)$-order, predict and compensate for the extracted signals, and reconstruct a new ADRCT that fully satisfies the separation principle in modern control theory. Finally, the numerical simulation and analysis results are given.

In this paper, the m-dimensional state space based on Euclidean straight space is used to derive the $m(m \geq 2)$-order synthesis function of discrete system by default, which will not be mentioned later.

In the latter part of this paper, the dynamic characteristics of system are mainly referred to when analyzing the zero-pole variation of transfer function, and when analyzing the phase delay and amplitude attenuation of extracted signals, the steady-state characteristics of system are mainly referred to.

## 2. Basic concepts and some combinatorial identities

The generating function is an important method in modern discrete mathematics. It is a bridge connecting discrete mathematics with continuous mathematics. It can effectively solve the combinatorial identity problems related to the combination numbers.

This section provides a brief introduction to the knowledge of generating function that will be used in the derivation of m-order synthesis function (Ke & Wei, 2010; Lint & Wilson, 2001; Wang, 2008), lists the combinatorial identities in theorem form to be used later, and then proves them in detail.

**Definition 2.1.** The function $sign(x)$ is defined by

$$sign(x) \triangleq \begin{cases} +1, x > 0 \\ 0, \quad x = 0 \\ -1, x < 0 \end{cases} \tag{1}$$

**Definition 2.2.** The saturation function is defined by

$$sat(x,\delta) \triangleq \begin{cases} sign(x), |x| > \delta \\ x/\delta, \quad |x| \leq \delta \end{cases} \tag{2}$$

**Definition 2.3.** The factorial of $k$ is defined by

$$k! \triangleq k \cdot (k-1) \cdots 2 \cdot 1 \tag{3}$$

**Definition 2.4.** $C_n^k$ is defined by

$$C_n^k \triangleq \frac{n!}{k!(n-k)!} \tag{4}$$

**Theorem 2.5.** If $1 \leq k < n$, then

$$C_k^0 = C_{k-1}^0 = 1 \tag{5}$$

$$C_k^n = 0 \tag{6}$$

$$C_n^k = C_{n-1}^k + C_{n-1}^{k-1} \tag{7}$$

**Definition 2.6.** Let $u \triangleq (u_i)_{i \geq 0} \triangleq (u_0, u_1, \cdots, u_n, \cdots)$ be an infinite series, then the formal power series (8) is the ordinary generating function of series $u$, referred as the generating function.

$$u(x) \triangleq \sum_{i \geq 0} u_i x^i \tag{8}$$

The "formal power series" in Definition 2.6 refers to not consider the convergence problem of series (8).

**Definition 2.7.** To say that two series are equal is to say that their coefficient sequences are same. For instance, the formal power series (8) is equal to series (9), if and only if (10) holds.

$$v(x) \triangleq \sum_{i\geq 0} v_i x^i \tag{9}$$

$$u_i = v_i \left(i \geq 0\right) \tag{10}$$

**Definition 2.8.** The operation that the number $a$ multiplies the series $u(x)$ is defined by

$$au(x) \triangleq a \cdot u(x) \triangleq \sum_{i\geq 0} a u_i x^i \tag{11}$$

The addition operation of series $u(x)$ and $v(x)$ is defined by

$$u(x)+v(x) \triangleq \sum_{i\geq 0}(u_i + v_i)x^i \tag{12}$$

The subtraction operation of series $u(x)$ and $v(x)$ is defined by

$$u(x)-v(x) \triangleq \sum_{i\geq 0}(u_i - v_i)x^i \tag{13}$$

The multiplication operation of series $u(x)$ and $v(x)$ is defined by

$$u(x)v(x) \triangleq u(x)\cdot v(x) \triangleq \sum_{k\geq 0}\left(\sum_{i+j=k} u_i v_j\right)x^k \tag{14}$$

**Lemma 2.9.** The generating function

$$\left(1+x\right)^n = \sum_{k=0}^{n} C_n^k x^k \tag{15}$$

**Lemma 2.10.** The generating function

$$\left(1-x\right)^n = \frac{1}{\left(1+x+x^2+\cdots\right)^n} = \sum_{k=0}^{n}(-1)^k C_n^k x^k \tag{16}$$

**Lemma 2.11.** The generating function

$$\frac{1}{\left(1-x\right)^n} = \left(1+x+x^2+\cdots\right)^n = \sum_{k=0}^{\infty} C_{n+k-1}^k x^k \tag{17}$$

**Lemma 2.12.** The generating function

$$\frac{1}{1+x} = 1-x+x^2-\cdots = \sum_{k=0}^{\infty}\left(-1\right)^k x^k \tag{18}$$

**Lemma 2.13.** The generating function

$$\frac{1}{\left(1+x\right)^n} = \left(1-x+x^2-\cdots\right)^n = \sum_{k=0}^{\infty}\left(-1\right)^k C_{n+k-1}^k x^k \tag{19}$$

**Theorem 2.14.** If $k \geq 0$, $m \geq 1$, $1 \leq i \leq m$, then there is identity

$$\sum_{j=1}^{k} C_{j+m-i-1}^{m-i} = \sum_{\mu=0}^{k-1} C_{\mu+m-i}^{m-i} = C_{m+k-i}^{m-i+1} \tag{20}$$

**Proof.** When $k = 0$

$$\sum_{j=1}^{k} C_{j+m-i-1}^{m-i} = 0 \tag{21}$$

$$C_{m+k-i}^{m-i+1} = C_{m-i}^{m-i+1} = 0 \tag{22}$$

$$\Rightarrow \sum_{j=1}^{k} C_{j+m-i-1}^{m-i} = C_{m+k-i}^{m-i+1} = 0 \tag{23}$$

When $k \geq 1$

$C_{j+m-i-1}^{m-i}$ is the coeff of item $x^{m-i}$ in generating function $\left(1+x\right)^{j+m-i-1}$, $\Rightarrow \sum_{j=1}^{k} C_{j+m-i-1}^{m-i}$ is the coeff of item $x^{m-i}$ in $\sum_{j=1}^{k}\left(1+x\right)^{j+m-i-1}$. While

$$\sum_{j=1}^{k}\left(1+x\right)^{j+m-i-1} = \left(1+x\right)^{m-i}\sum_{j=0}^{k-1}\left(1+x\right)^j$$

$$= \left(1+x\right)^{m-i}\left[\frac{\left(1+x\right)^k - 1}{1+x-1}\right]$$

$$= \frac{1}{x}\left[\left(1+x\right)^{m+k-i} - \left(1+x\right)^{m-i}\right] \tag{24}$$

$$\Rightarrow \sum_{j=1}^{k} C_{j+m-i-1}^{m-i} = C_{m+k-i}^{m-i+1} \tag{25}$$

□

**Theorem 2.15.** If $k, m \geq 1$, then there is identity

$$\sum_{i=0}^{m-1}\left(-1\right)^i C_{m-1}^i C_{k+m-i-1}^{m-i} = C_k^m \tag{26}$$

**Proof.** $\left(-1\right)^i C_{m-1}^i$ is the coeff of item $x^i$ in $f_1(x) = \left(1-x\right)^{m-1}$, $C_{k+m-i-1}^{m-i}$ is the coeff of item $x^{m-i}$ in $f_2(x) = \left(1+x+x^2+\cdots\right)^k$, $\Rightarrow \sum_{i=0}^{m}\left(-1\right)^i C_{m-1}^i C_{k+m-i-1}^{m-i}$ is the coeff of item $x^i x^{m-i} = x^m$ in $f(x) = f_1(x) f_2(x)$.

When $k \geq m \geq 1$

$$f(x) = \left(1-x\right)^{m-1}\left(1+x+x^2+\cdots\right)^k$$

$$= \left(1+x+x^2+\cdots\right)^{k-(m-1)} \tag{27}$$

$$\Rightarrow \sum_{i=0}^{m}\left(-1\right)^i C_{m-1}^i C_{k+m-i-1}^{m-i} = C_{k-(m-1)+m-1}^m = C_k^m \tag{28}$$

When $1 \le k \le m-1$

$$f(x) = (1-x)^{(m-1)-k} \tag{29}$$

$$\Rightarrow \sum_{i=0}^{m} (-1)^i C_{m-1}^i C_{k+m-i-1}^{m-i} = 0 \tag{30}$$

Notice equation (6), when $m \ge 1$, equation (26) holds. □

**Theorem 2.16.** If $k \ge m-1, m \ge 2$, then there are identities

$$\sum_{i=0}^{m-1} (-1)^i C_k^i C_{k+m-i-1}^{m-i} = (-1)^{m-1} C_k^m \tag{31}$$

$$\sum_{i=0}^{m-1} (-1)^i C_k^i C_{k+m-i-2}^{m-i} = (-1)^{m-1} C_k^m \tag{32}$$

**Proof.** $(-1)^i C_k^i$ is the coeff of item $x^i$ in $f_1(x) = (1-x)^k$, $C_{k+m-i-1}^{m-i}$ is the coeff of item $x^{m-i}$ in $f_2(x) = (1+x+x^2+\cdots)^k$, $\Rightarrow \sum_{i=0}^{m} (-1)^i C_k^i C_{k+m-i-1}^{m-i}$ is the coeff of item $x^m$ in $f(x) = f_1(x) f_2(x)$.

$$f(x) = (1-x)^k (1+x+x^2+\cdots)^k = 1 \tag{33}$$

$$\Rightarrow \sum_{i=0}^{m} (-1)^i C_k^i C_{k+m-i-1}^{m-i} = 0 \tag{34}$$

$$\Rightarrow \sum_{i=0}^{m-1} (-1)^i C_k^i C_{k+m-i-1}^{m-i} + (-1)^m C_k^m C_{k-1}^0 = 0 \tag{35}$$

$$\Rightarrow \sum_{i=0}^{m-1} (-1)^i C_k^i C_{k+m-i-1}^{m-i} = (-1)^{m-1} C_k^m \tag{36}$$

Similarly, we can prove (32). □

**Theorem 2.17.** If $k, m \ge 1$, then there is identity

$$\sum_{i=0}^{m-1} (-1)^i C_k^i = (-1)^{m-1} C_{k-1}^{m-1} \tag{37}$$

**Proof.** We first prove $\sum_{i=0}^{m} (-1)^{m-i} C_k^i = C_{k-1}^m$.

$(-1)^{m-i}$ is the coeff of item $x^{m-i}$ in $\frac{1}{1+x}$, $C_k^i$ is the coeff of item $x^i$ in $(1+x)^k$, $\Rightarrow \sum_{i=0}^{m} (-1)^{m-i} C_k^i$ is the coeff of item $x^m$ in $\frac{1}{1+x}(1+x)^k = (1+x)^{k-1}$.

$$\Rightarrow \sum_{i=0}^{m} (-1)^{m-i} C_k^i = C_{k-1}^m \tag{38}$$

$$\Rightarrow \sum_{i=0}^{m-1} (-1)^{m-i} C_k^i + (-1)^0 C_k^m = C_{k-1}^m \tag{39}$$

Notice equation (7), then there is

$$C_k^m = C_{k-1}^m + C_{k-1}^{m-1} \tag{40}$$

$$\Rightarrow \sum_{i=0}^{m-1} (-1)^{m-i} C_k^i = -C_{k-1}^{m-1} \tag{41}$$

$$\Rightarrow \sum_{i=0}^{m-1} (-1)^i C_k^i = (-1)^m \sum_{i=0}^{m-1} (-1)^{m-i} C_k^i = (-1)^{m-1} C_{k-1}^{m-1} \tag{42}$$

□

**Theorem 2.18.** If $k, m \ge 1$, then there is identity

$$\sum_{i=0}^{m-1} (-1)^i C_m^i C_{k+m-i-1}^{m-i} = C_{k-1}^m + (-1)^{m-1} \tag{43}$$

**Proof.** $(-1)^i C_m^i$ is the coeff of item $x^i$ in $f_1(x) = (1-x)^m$, $C_{k+m-i-1}^{m-i}$ is the coeff of item $x^{m-i}$ in $f_2(x) = (1+x+x^2+\cdots)^k$, $\Rightarrow \sum_{i=0}^{m} (-1)^i C_m^i C_{k+m-i-1}^{m-i}$ is the coeff of item $x^m$ in $f(x) = f_1(x) f_2(x)$.

When $k \ge m+1$

$$f(x) = (1-x)^m (1+x+x^2+\cdots)^k = (1+x+x^2+\cdots)^{k-m} \tag{44}$$

$$\Rightarrow \sum_{i=0}^{m} (-1)^i C_m^i C_{k+m-i-1}^{m-i} = C_{k-m+m-1}^m = C_{k-1}^m \tag{45}$$

When $1 \le k \le m$

$$f(x) = (1-x)^{m-k}$$

$$\Rightarrow \sum_{i=0}^{m} (-1)^i C_m^i C_{k+m-i-1}^{m-i} = 0 \tag{46}$$

Notice equation (6), when $m \ge 1$, equation (45) holds.

$$\sum_{i=0}^{m-1} (-1)^i C_m^i C_{k+m-i-1}^{m-i} + (-1)^m C_m^m C_{k-1}^0 = C_{k-1}^m \tag{47}$$

$$\Rightarrow \sum_{i=0}^{m-1} (-1)^i C_m^i C_{k+m-i-1}^{m-i} = C_{k-1}^m + (-1)^{m-1} \tag{48}$$

□

**Theorem 2.19.** If $m \ge 1$, then there is identity

$$\sum_{i=0}^{m-1} (-1)^i C_m^i = (-1)^{m-1} \tag{49}$$

**Proof.**

$$(1-x)^m = \sum_{i=0}^{m} (-1)^i C_m^i x^i \tag{50}$$

Let $x = 1$, then obtain

$$\sum_{i=0}^{m}(-1)^i C_m^i = 0 \tag{51}$$

$$\Rightarrow \sum_{i=0}^{m-1}(-1)^i C_m^i + (-1)^m C_m^m = 0 \tag{52}$$

$$\Rightarrow \sum_{i=0}^{m-1}(-1)^i C_m^i = (-1)^{m-1} \tag{53}$$

□

**Theorem 2.20.** If $m \geq 1$, $0 \leq \nu \leq m-1$, $0 \leq k \leq m-1-\nu$, then there is identity

$$\sum_{i=0}^{m-(\nu+1)}(-1)^i C_{m-(\nu+1)}^i C_{k+m-(i+\nu)-1}^{m-(i+\nu)} = 0 \tag{54}$$

**Proof.** When $k=0$

$C_{k+m-(i+\nu)-1}^{m-(i+\nu)} = C_{m-(i+\nu)-1}^{m-(i+\nu)} = 0$, equation (54) holds.

When $1 \leq k \leq m-1-\nu$

$(-1)^i C_{m-(\nu+1)}^i$ is the coeff of item $x^i$ in $f_1(x) = (1-x)^{m-(\nu+1)}$, $C_{k+m-(i+\nu)-1}^{m-(i+\nu)}$ is the coeff of item $x^{m-(i+\nu)}$ in $f_2(x) = (1+x+x^2+\cdots)^k$,

$\Rightarrow \sum_{i=0}^{m-\nu}(-1)^i C_{m-(\nu+1)}^i C_{k+m-(i+\nu)-1}^{m-(i+\nu)}$ is the coeff of item $x^{m-\nu}$ in $f(x) = f_1(x) f_2(x)$.

$$f(x) = (1-x)^{m-(\nu+1)-k} \tag{55}$$

$$\sum_{i=0}^{m-\nu}(-1)^i C_{m-(\nu+1)}^i C_{k+m-(i+\nu)-1}^{m-(i+\nu)} = (-1)^{m-\nu} C_{m-(\nu+1)-k}^{m-\nu} = 0 \tag{56}$$

Notice equation (6), when $m \geq 1$, equation (54) holds. □

## 3. Second-order synthesis function

In order to understand the derivation source and control ideas of m-order synthesis function, this section first introduces the derivation process and control ideas of 2nd-order synthesis function $fsun()$. In order to maintain the integrity described in this section, this section cites the discussion in (Sun & Sun, 2010).

The 2nd-order discrete system satisfies the following equations.

$$\begin{cases} x_1(k+1) = x_1(k) + hx_2(k) \\ x_2(k+1) = x_2(k) + hu(k) \end{cases} \tag{57}$$

Where, $h$ is the sampling time (or control step), $u$ is the control input of system (57), $|u(k)| \leq r$, $r$ is the maximum control input.

Now we need to solve the time optimal control synthesis function of system (57) based on the state variables (i.e., the expression of time optimal control input $u$ based on the state variables) in the following.

### 3.1. Time optimal trajectory

Set the initial value is $X(0) \triangleq [x_1(0), x_2(0)]^T$. After $k$ steps, the solution of system (57) is

$$X(k) \triangleq [x_1(k), x_2(k)]^T =$$

$$\begin{bmatrix} 1 & kh \\ 0 & 1 \end{bmatrix}\begin{bmatrix} x_1(0) \\ x_2(0) \end{bmatrix} + \begin{bmatrix} 1 & (k-1)h \\ 0 & 1 \end{bmatrix}\begin{bmatrix} 0 \\ h \end{bmatrix} u(0) +$$

$$\cdots + \begin{bmatrix} 1 & h \\ 0 & 1 \end{bmatrix}\begin{bmatrix} 0 \\ h \end{bmatrix} u(k-2) + \begin{bmatrix} 0 \\ h \end{bmatrix} u(k-1) \tag{58}$$

Assume that system (57) reaches the steady-state (i.e., origin) after $k$ steps control, that is

$$\begin{bmatrix} x_1(k) \\ x_2(k) \end{bmatrix} = \begin{bmatrix} 0 \\ 0 \end{bmatrix} \tag{59}$$

The all initial values $X(0)$ of system that can reach the origin within $k$ steps (the set of all these points is called the $k$ isochronous time region in 2nd-dimensional state space (Qian & Song, 2011), denoted as $G_2(k)$ ) are expressed as

$$X(0) = \begin{bmatrix} x_1(0) \\ x_2(0) \end{bmatrix} = \begin{bmatrix} kh^2 \\ -h \end{bmatrix} u(k-1) + \begin{bmatrix} (k-1)h^2 \\ -h \end{bmatrix} u(k-2) +$$

$$\cdots + \begin{bmatrix} 2h^2 \\ -h \end{bmatrix} u(1) + \begin{bmatrix} h^2 \\ -h \end{bmatrix} u(0) \tag{60}$$

The point $a_{-k}\ (k \geq 1)$ is defined as the initial point that can reach the origin along the polyline $a_{-k}$, $a_{-(k-1)}$, $\cdots$, $a_{-1}$, 0 in $k$ steps under all control inputs $u$ taking $-r$. Similarly, the point $a_{+k}\ (k \geq 1)$ is defined as the initial point that can reach the origin along the polyline $a_{+k}$, $a_{+(k-1)}$, $\cdots$, $a_{+1}$, 0 in $k$ steps under all controls $u$ taking $+r$.

The point $b_{+k}\ (k \geq 2)$ is defined as the initial point that can first reach the point $a_{-(k-1)}$ under the control $u$ taking $+r$, and can then reach the origin along the polyline $a_{-(k-1)}$, $a_{-(k-2)}$, $\cdots$, $a_{-1}$, 0 in $(k-1)$ steps under all other controls $u$ taking $-r$. Similarly, the point $b_{-k}\ (k \geq 2)$ is defined as the initial point that can first reach the point $a_{+(k-1)}$ under the control $u$ taking $-r$, and can then reach the origin along the polyline $a_{+(k-1)}$, $a_{+(k-2)}$, $\cdots$, $a_{+1}$, 0 in $(k-1)$ steps under all other controls $u$ taking $+r$.

The point $a_0$ is defined as the origin.

Thus, the corresponding coordinates of points ($a_{-k}$, $a_{+k}$, $b_{+k}$, and $b_{-k}$) can be obtained as follows.

The coordinate of point $a_{-k}$ is

$$\begin{bmatrix} x_1 \\ x_2 \end{bmatrix} = \begin{bmatrix} -\frac{1}{2}k(k+1)h^2r \\ khr \end{bmatrix} \tag{61}$$

The coordinate of point $a_{+k}$ is

$$\begin{bmatrix} x_1 \\ x_2 \end{bmatrix} = \begin{bmatrix} \frac{1}{2}k(k+1)h^2r \\ -khr \end{bmatrix} \tag{62}$$

The coordinate of point $b_{+k}$ is

$$\begin{bmatrix} x_1 \\ x_2 \end{bmatrix} = \begin{bmatrix} -\frac{1}{2}k(k+1)h^2r + 2h^2r \\ (k-2)hr \end{bmatrix} \tag{63}$$

The coordinate of point $b_{-k}$ is

$$\begin{bmatrix} x_1 \\ x_2 \end{bmatrix} = \begin{bmatrix} \frac{1}{2}k(k+1)h^2r - 2h^2r \\ -(k-2)hr \end{bmatrix} \tag{64}$$

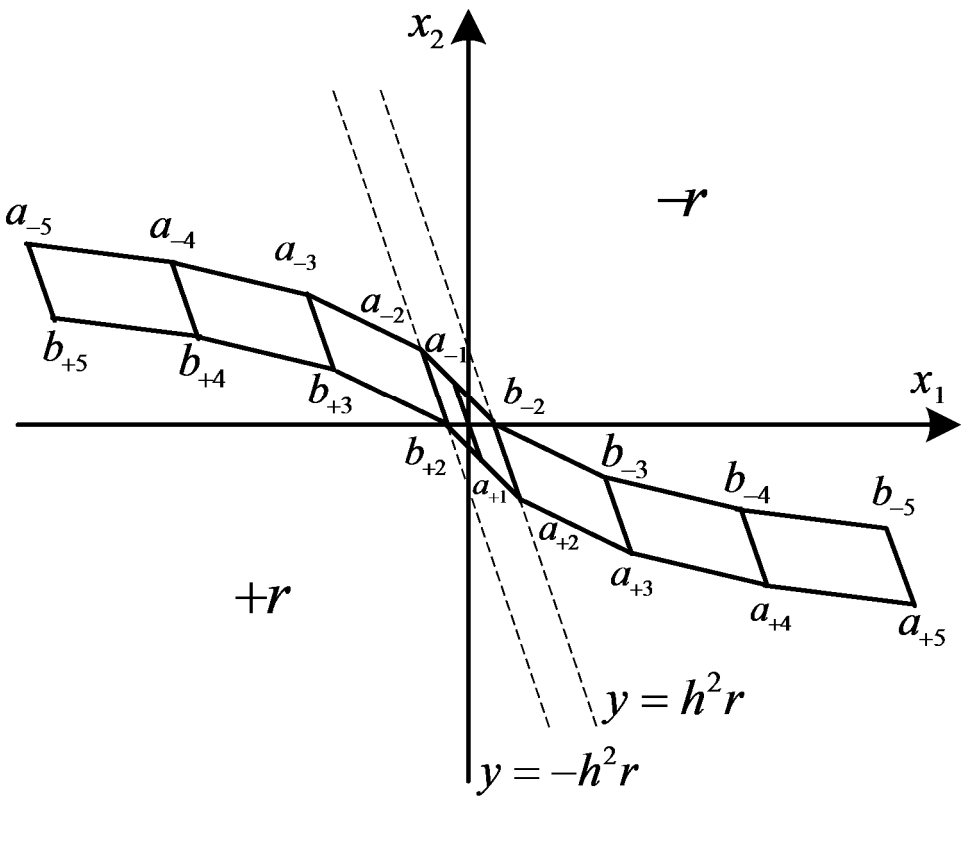

Fig.1. Time optimal trajectory of 2nd-order discrete system

The distribution of points ($a_{-k}$, $a_{+k}$, $b_{+k}$, and $b_{-k}$) are shown in Fig.1. In Fig.1, the polyline formed by $a_{+k}$, $a_{+(k-1)}$, $\cdots$, $a_{+1}$, 0, $a_{-1}$, $\cdots$, $a_{-(k-1)}$, $a_{-k}$ is denoted as $G_{best}$, and obviously $G_{best}$ is a kind of time optimal trajectory to reach the origin. The polyline formed by $a_{+k}$, $a_{+(k-1)}$, $\cdots$, $a_{+1}$, $b_{+2}$, $\cdots$, $b_{+(k-1)}$, $b_{+k}$ is denoted as $G_+$, and the polyline formed by $a_{-k}$, $a_{-(k-1)}$, $\cdots$, $a_{-1}$, $b_{-2}$, $\cdots$, $b_{-(k-1)}$, $b_{-k}$ is denoted as $G_-$.

Obviously, the control input $u$ takes the value within $[-r,+r]$ if the state variable $X$ is located within the region enclosed by the polyline $G_-$ and polyline $G_+$, the control $u$ takes $+r$ if the variable $X$ is located on the polyline $G_+$ or below (outside) $G_+$, and $u$ takes $-r$ if the variable $X$ is located on the polyline $G_-$ or above (outside) $G_-$.

**Definition 3.1.** The functions are defined as follows.

$$y \triangleq x_1 + hx_2 \tag{65}$$

$$z(k) \triangleq x_1 + khx_2 \tag{66}$$

$$s \triangleq sign(y) \tag{67}$$

In Fig.1, the line equation (denoted as $N_2(+k)$, where $y > 0$) through the line segment $a_{+k}b_{-k}$ $(k \ge 2)$ is

$$y = x_1 + hx_2 = \frac{1}{2}k(k-1)h^2r \tag{68}$$

The line equation (denoted as $N_2(-k)$, where $y < 0$) through the line segment $a_{-k}b_{+k}$ $(k \ge 2)$ is

$$y = x_1 + hx_2 = -\frac{1}{2}k(k-1)h^2r \tag{69}$$

After $a_k$ uniformly represents $a_{+k}$ and $a_{-k}$, $b_k$ uniformly represents $b_{-k}$ and $b_{+k}$, then the line equations (68) and (69) (uniformly denoted as $N_2(k)$) are uniformly expressed as

$$y = x_1 + hx_2 = \frac{1}{2}k(k-1)h^2rs \tag{70}$$

The line segment $a_ka_{k-1}$ (uniformly denoted as $M_2(k)$: the line segment $a_{-k}a_{-(k-1)}$ is denoted as $M_2(+k)$, where any point on the line segment satisfies $z(k) > 0$, $s < 0$; and the line segment $a_{+k}a_{+(k-1)}$ is denoted as $M_2(-k)$, where any point on the line segment satisfies $z(k) < 0$, $s > 0$.) satisfies the line equation as follows.

$$z(k) = x_1 + khx_2 = -\frac{1}{2}k(k-1)h^2rs \tag{71}$$

Similarly, the line segment $a_{k-1}a_{k-2}$ (denoted as $M_2(k-1)$) satisfies the line equation as follows.

$$z(k-1) = x_1' + (k-1)hx_2' = -\frac{1}{2}(k-1)(k-2)h^2rs \tag{72}$$

### *3.2. The nonlinear region of $|y| > h^2 r$*

Suppose that the state $X \triangleq [x_1, x_2]^T$ is located between the parallel lines $N_2(k)$ and $N_2(k-1)$, satisfying the following equation.

$$y = x_1 + hx_2 = \frac{1}{2}k'(k'-1)h^2 rs \tag{73}$$

Where $k$ is a positive integer and satisfies

$$\begin{cases} k-1 < k' \le k \\ k \ge 2 \end{cases} \tag{74}$$

Obtain from (73)

$$k'(k'-1) = \frac{2y}{h^2 rs} = \frac{2|y|}{h^2 r} \tag{75}$$

Therefore, the value $k$ meets

$$\begin{cases} k'(k'-1) = \dfrac{2|y|}{h^2 r} \\ k-1 < k' \le k \\ k \ge 2 \end{cases} \tag{76}$$

The state $X$ is located between the parallel lines $N_2(k)$ and $N_2(k-1)$, and the control target is that the state $X$ reaches the point (set to $X' \triangleq [x_1', x_2']^T$ ) on the line segment $a_{k-1}a_{k-2}$ after one-step control. The new state $X'$ meets

$$\begin{cases} x_1' = x_1 + hx_2 \\ x_2' = x_2 + hu \end{cases} \tag{77}$$

Substitute (77) into (72) and obtain

$$x_1 + khx_2 + (k-1)h^2 u = -\frac{1}{2}(k-1)(k-2)h^2 rs \tag{78}$$

$$\Rightarrow u = \left(1 - \frac{k}{2}\right) rs - \frac{x_1 + khx_2}{(k-1)h^2}, \quad |y| > h^2 r \tag{79}$$

**Definition 3.2.** $k \ge 2$ , $|y| > h^2 r$ , the function $a(x_1, x_2, r, h)$ is defined as follows.

$$a(x_1, x_2, r, h) \triangleq -\left(1 - \frac{k}{2}\right) rs + \frac{x_1 + khx_2}{(k-1)h^2} \tag{80}$$

Since $|u| \le r$ and $|y| > h^2 r$ , then the control input (i.e., the 2nd-order synthesis function of nonlinear region) is

$$u = -r \cdot sat\left[a(x_1, x_2, r, h), r\right] \tag{81}$$

### *3.3. The linear region of $|y| \le h^2 r$*

If the state $X$ satisfies $|y| \le h^2 r$ , that is, $X$ is located between the parallel lines $N_2(+2)$ and $N_2(-2)$, and the control target is that the state $X$ reaches the point (set to $X' \triangleq [x_1', x_2']^T$ ) on the line segment $a_{+1}a_{-1}$ and then reaches the origin along the line segment $a_1 a_0$ .

The line segment $a_{+1}a_{-1}$ (denoted as $N_2(1)$ or $M_2(1)$ ) satisfies the following line equation.

$$y = x_1' + hx_2' = 0 \tag{82}$$

Substitute (77) into (82) and obtain

$$x_1 + 2hx_2 + h^2 u = 0 \tag{83}$$

$$\Rightarrow u = -\frac{x_1 + 2hx_2}{h^2}, \quad |y| \le h^2 r \tag{84}$$

For $|u| \le r$ and $|y| \le h^2 r$ , then the control input (i.e., the 2nd-order synthesis function of linear region) is

$$u = -r \cdot sat\left(\frac{x_1 + 2hx_2}{h^2}, r\right) \tag{85}$$

### *3.4. Second-order synthesis function*

**Definition 3.3.** In both cases of $|y| > h^2 r$ and $|y| \le h^2 r$ , the function $a(x_1, x_2, r, h)$ is uniformly defined as follows.

$$a(x_1, x_2, r, h) \triangleq \begin{cases} \dfrac{x_1 + 2hx_2}{h^2}, & |y| \le h^2 r \\ -\left(1 - \dfrac{k}{2}\right) rs + \dfrac{x_1 + khx_2}{(k-1)h^2}, & |y| > h^2 r \end{cases} \tag{86}$$

Where $k$ is a positive integer and meets (76).

**Theorem 3.4.** The 2nd-order time optimal control input of discrete system (57) (denoted as $fsun(x_1, x_2, r, h)$ ) is

$$u = fsun(x_1, x_2, r, h) \triangleq -r \cdot sat\left[a(x_1, x_2, r, h), r\right] \tag{87}$$

(87) is the 2nd-order synthesis function of discrete system (57) with the origin as the end point. Obviously, this is a nonlinear function.

**Theorem 3.5.** If the control input limit $|u| \le r$ is removed in the 2nd-order discrete system (57), i.e., $r \to +\infty$ , then we can get its 2nd-order linear synthesis function with the origin as the end point as follows.

$$u = -\frac{x_1 + 2hx_2}{h^2} \tag{88}$$

## 4. m-order synthesis function

The $m(m \ge 2)$ -order discrete system satisfies

$$\begin{cases} x_1(k+1) = x_1(k) + hx_2(k) \\ x_2(k+1) = x_2(k) + hx_3(k) \\ \vdots \\ x_m(k+1) = x_m(k) + hu(k) \end{cases} \quad (89)$$

Where, $h$ is the sampling time (or control step), $u$ is the control input of system (89), $|u(k)| \le r$, $r$ is the maximum control input.

Now we need to solve the time optimal control synthesis function of system (89) based on the state variables (i.e., the expression of time optimal control input $u$ based on the state variables) in the following.

The system (89) is written in matrix form as follows.

$$X(k+1) = AX(k) + B_0 u(k) \quad (90)$$

Where, $X(k) \triangleq [x_1(k), x_2(k), \cdots, x_m(k)]^T$ (91)

$$A = \begin{bmatrix} 1 & h & 0 & & 0 \\ 0 & 1 & h & \cdots & 0 \\ 0 & 0 & 1 & & 0 \\ & \vdots & & \ddots & \vdots \\ 0 & 0 & 0 & & h \\ 0 & 0 & 0 & \cdots & 1 \end{bmatrix}, \ B_0 = \begin{bmatrix} 0 \\ 0 \\ 0 \\ \vdots \\ 0 \\ h \end{bmatrix} \quad (92)$$

Set $A_k \triangleq A^k$, $B_{k-1} \triangleq A^{k-1} B_0$ (93)

### *4.1. Solution of equations*

Set the initial value is $X(0) \triangleq [x_1(0), x_2(0), \cdots, x_m(0)]^T$. After $k$ steps, the solution of system (90) is

$$X(k) = A_k X(0) + \sum_{\mu=0}^{k-1} B_{(k-1-\mu)} u(\mu) \quad (94)$$

Assume that system (90) reaches the origin after $k$ steps control, i.e.,

$$X(k) = [0, 0, \cdots, 0]^T \quad (95)$$

The all initial values $X(0)$ of system that can reach the origin within $k$ steps (the set of all these points is called the $k$ isochronous time region in m-dimensional state space (Qian & Song, 2011), denoted as $G_m(k)$) are expressed as

$$X(0) = -A_k^{-1} \sum_{\mu=0}^{k-1} B_{(k-1-\mu)} u(\mu) \quad (96)$$

Set $C_k^{j-i} \triangleq 0$, when $j < i$ (97)

The corresponding matrix expression that we can get is

$$A_k(i,j) = h^{j-i} \times \left( coeff \ of \ item \ x^j \ in \ x^i (1+x)^k \right)$$

$$= C_k^{j-i} h^{j-i}, \ 1 \le i, j \le m \quad (98)$$

$$B_{(k-1-\mu)} = A^{k-1-\mu} B_0 \quad (99)$$

$$B_{(k-1-\mu)}(i) = hA_{k-1-\mu}(i,m) = C_{k-1-\mu}^{m-i} h^{m-i+1}, \ 1 \le i \le m \quad (100)$$

**Theorem 4.1.** $k \ge 1$, $1 \le i, j \le m$, the inverse matrix of $A_k(i,j) = C_k^{j-i} h^{j-i}$ is

$$A_k^{-1}(i,j) = h^{j-i} \left( coeff \ of \ item \ x^j \ in \ \frac{x^i}{(1+x)^k} \right)$$

$$= (-1)^{j-i} C_{k+j-i-1}^{j-i} h^{j-i}, \ 1 \le i, j \le m \quad (101)$$

**Proof.**

$$A_k(i,\nu) = C_k^{\nu-i} h^{\nu-i}, \ 1 \le i, \nu \le m \quad (102)$$

Set $A'(\nu, j) = h^{j-\nu} \left( coeff \ of \ item \ x^j \ in \ \dfrac{x^\nu}{(1+x)^k} \right)$

$$= (-1)^{j-\nu} C_{k+j-\nu-1}^{j-\nu} h^{j-\nu}, \ 1 \le \nu, j \le m \quad (103)$$

$$\Rightarrow \sum_{\nu=1}^{m} A_k(i,\nu) A'(\nu,j) = h^{j-i} \sum_{\nu=1}^{m} (-1)^{j-\nu} C_k^{\nu-i} C_{k+j-\nu-1}^{j-\nu} \quad (104)$$

When $i = j$, since $C_k^{\nu-i}$ and $C_{k+j-\nu-1}^{j-\nu}$ in (104) have exact definitions at the same time, $\Rightarrow \nu - i \ge 0$, $j - \nu \ge 0$, $\Rightarrow \nu = i = j \Rightarrow$ (104) is equal to 1;

When $i > j$, since $C_k^{\nu-i}$ and $C_{k+j-\nu-1}^{j-\nu}$ in (104) have exact definitions at the same time, $\Rightarrow \nu - i \ge 0$, $j - \nu \ge 0$ $\Rightarrow i \le \nu \le j$, it is contradictory with the hypothesis $i > j \Rightarrow$ (104) is equal to 0;

When $i < j$, $C_k^{\nu-i}$ and $C_{k+j-\nu-1}^{j-\nu}$ in (104) have exact definitions at the same time only while $i \le \nu \le j$,

$$\sum_{\nu=1}^{m} A_k(i,\nu) A'(\nu,j)$$

$$= h^{j-i} \sum_{\nu=i}^{j} (-1)^{j-\nu} \frac{k!}{(\nu-i)!(k+i-\nu)!} \frac{(k+j-\nu-1)!}{(j-\nu)!(k-1)!}$$

$$= \frac{kh^{j-i}}{j-i} \sum_{\nu=i}^{j} (-1)^{j-\nu} \frac{(k+j-\nu-1)!}{(j-i-1)!(k+i-\nu)!} \frac{(j-i)!}{(j-\nu)!(\nu-i)!}$$

$$= \frac{kh^{j-i}}{j-i} \sum_{\nu=i}^{j} (-1)^{j-\nu} C_{j-i+(k+i-\nu)-1}^{k+i-\nu} C_{j-i}^{\nu-i}$$

$$= \frac{kh^{j-i}}{j-i} (-1)^{k+j-i} \sum_{\nu=i}^{j} (-1)^{k+i-\nu} C_{j-i+(k+i-\nu)-1}^{k+i-\nu} C_{j-i}^{\nu-i} \quad (105)$$

In (105), $\left(-1\right)^{k+i-\nu} C_{j-i+(k+i-\nu)-1}^{k+i-\nu}$ is the coeff of item $x^{k+i-\nu}$ in $\frac{1}{\left(1+x\right)^{j-i}}$, $C_{j-i}^{\nu-i}$ is the coeff of item $x^{\nu-i}$ in $\left(1+x\right)^{j-i}$, and therefore $\sum_{\nu=i}^{j}\left(-1\right)^{k+i-\nu} C_{j-i+(k+i-\nu)-1}^{k+i-\nu} C_{j-i}^{\nu-i}$ is the coeff of item $x^{k}$ in $\frac{1}{\left(1+x\right)^{j-i}}\left(1+x\right)^{j-i}=1$. For $k\ge 1$, so (105) is equal to 0.

Overall, there is

$$\sum_{\nu=1}^{m} A_k\left(i,\nu\right)A'\left(\nu,j\right)=\delta_{ij} \tag{106}$$

□

**Theorem 4.2.** $k\ge 1$, $0\le \mu\le k-1$, $1\le i,j\le m$

$$\sum_{j=1}^{m} A_k^{-1}\left(i,j\right)B_{(k-1-\mu)}\left(j\right)=\sum_{j=1}^{m}\left(-1\right)^{j-i} C_{k+j-i-1}^{j-i}h^{j-i}C_{k-1-\mu}^{m-j}h^{m-j+1}$$

$$=h^{m-i+1}\sum_{j=1}^{m}\left(-1\right)^{j-i} C_{k+j-i-1}^{j-i}C_{k-1-\mu}^{m-j}$$

$$=\left(-1\right)^{m-i} C_{m-i+\mu}^{m-i}h^{m-i+1} \tag{107}$$

**Proof.**

In (107), $\left(-1\right)^{j-i} C_{k+j-i-1}^{j-i}$ is the coeff of item $x^{j-i}$ in $\frac{1}{\left(1+x\right)^{k}}$, $C_{k-1-\mu}^{m-j}$ is the coeff of item $x^{m-j}$ in $\left(1+x\right)^{k-1-\mu}$, and $\sum_{j=1}^{m}\left(-1\right)^{j-i} C_{k+j-i-1}^{j-i}C_{k-1-\mu}^{m-j}$ is the coeff of item $x^{m-i}$ in $\frac{1}{\left(1+x\right)^{k}}\left(1+x\right)^{k-1-\mu}=\frac{1}{\left(1+x\right)^{\mu+1}}$, therefore

$$\sum_{j=1}^{m} A_k^{-1}\left(i,j\right)B_{(k-1-\mu)}\left(j\right)=\left(-1\right)^{m-i} C_{m-i+\mu}^{m-i}h^{m-i+1} \tag{108}$$

□

(96) can be written as follows.

$$x_i\left(0\right)=\left(-1\right)^{m-i+1}\sum_{\mu=0}^{k-1}C_{m-i+\mu}^{m-i}h^{m-i+1}u\left(\mu\right),\ \ 1\le i\le m \tag{109}$$

The points $a_{-k}\left(k\ge 1\right)$, $a_{+k}\left(k\ge 1\right)$, $b_{-k}\left(k\ge 2\right)$ and $b_{+k}\left(k\ge 2\right)$ are defined as section 3. The point $a_0$ is the origin.

By use of Theorem 2.14, the coordinate of point $a_{-k}\left(k\ge 1\right)$ can be obtained as follows.

$$x_i=\left(-1\right)^{m-i}\sum_{\mu=0}^{k-1}C_{m-i+\mu}^{m-i}h^{m-i+1}r$$

$$=\left(-1\right)^{m-i} C_{m+k-i}^{m-i+1}h^{m-i+1}r,\ \ 1\le i\le m \tag{110}$$

The coordinate of point $a_{+k}\left(k\ge 1\right)$ is

$$x_i=\left(-1\right)^{m-i+1}\sum_{\mu=0}^{k-1}C_{m-i+\mu}^{m-i}h^{m-i+1}r$$

$$=\left(-1\right)^{m-i+1} C_{m+k-i}^{m-i+1}h^{m-i+1}r,\ \ 1\le i\le m \tag{111}$$

The coordinate of point $b_{-k}\left(k\ge 2\right)$ is

$$x_i=\left(-1\right)^{m-i+1}\left(\sum_{\mu=0}^{k-1}C_{m-i+\mu}^{m-i}-2\right)h^{m-i+1}r$$

$$=\left(-1\right)^{m-i+1}\left(C_{m+k-i}^{m-i+1}-2\right)h^{m-i+1}r,\ \ 1\le i\le m \tag{112}$$

The coordinate of point $b_{+k}\left(k\ge 2\right)$ is

$$x_i=\left(-1\right)^{m-i}\left(\sum_{\mu=0}^{k-1}C_{m-i+\mu}^{m-i}-2\right)h^{m-i+1}r$$

$$=\left(-1\right)^{m-i}\left(C_{m+k-i}^{m-i+1}-2\right)h^{m-i+1}r,\ \ 1\le i\le m \tag{113}$$

*4.2. (m-1)-dimensional hyperplane $N_m\left(k\right)$*

**Definition 4.3.** The functions are defined as follows.

$$y\triangleq y\left(X\right)\triangleq\sum_{i=0}^{m-1}C_{m-1}^{i}h^{i}x_{i+1} \tag{114}$$

$$s\triangleq sign\left(y\right) \tag{115}$$

Unless we explicitly state that $s$ represents Laplace operator in the following analysis, (114) and (115) are used elsewhere to define the value of $s$.

**Definition 4.4.** The $\left(m-1\right)$-dimensional hyperplane $N_m\left(k\right)$ is defined by

$$y=C_k^{m}h^{m}rs \tag{116}$$

When $s>0$, $y>0$, $k\ge m$, $N_m\left(k\right)$ is denoted as $N_m\left(+k\right)$; when $s<0$, $y<0$, $k\ge m$, $N_m\left(k\right)$ is denoted as $N_m\left(-k\right)$; when $s=0$, $y=0$, $k\le m-1$, $N_m\left(k\right)$, is no longer distinguished between the positive and negative hyperplane, $N_m\left(m-1\right)$, $N_m\left(m-2\right)$, $N_m\left(1\right)$ are the same hyperplane and then are uniformly denoted as $N_m(m-1)$.

**Theorem 4.5.** 1) The points $a_{+k}\left(k\ge 1\right)$ and $b_{-k}\left(k\ge 2\right)$ are on the following $\left(m-1\right)$-dimensional hyperplane.

$$y=\left(-1\right)^{m} C_k^{m}h^{m}r \tag{117}$$

2) The points $a_{-k}\left(k\ge 1\right)$ and $b_{+k}\left(k\ge 2\right)$ are on the following $\left(m-1\right)$-dimensional hyperplane.

$$y=(-1)^{m-1}C_k^m h^m r \tag{118}$$

3) $$s=\begin{cases}(-1)^m, & a_{+k},b_{-k}\\(-1)^{m-1}, & a_{-k},b_{+k}\end{cases} \tag{119}$$

4) The points $a_{+k}(k\ge1)$ and $a_{-k}(k\ge1)$ are uniformly denoted as $a_k(k\ge1)$; $b_{-k}(k\ge2)$ and $b_{+k}(k\ge2)$ are uniformly denoted as $b_k(k\ge2)$. Then the points $a_k(k\ge1)$ and $b_k(k\ge2)$ are on the $(m-1)$-dimensional hyperplane $N_m(k)$.

**Proof.** 1) For the point $a_{+k}(k\ge1)$, substitute (111) into (114), notice Theorem 2.15 and obtain

$$y=\sum_{i=0}^{m-1}C_{m-1}^i h^i x_{i+1}=\sum_{i=0}^{m-1}C_{m-1}^i h^i(-1)^{m-i}C_{m+k-i-1}^{m-i}h^{m-i}r$$

$$=(-1)^m C_k^m h^m r \tag{120}$$

For the point $b_{-k}(k\ge2)$, substitute (112) into (114), notice Theorem 2.15 and also obtain

$$y=\sum_{i=0}^{m-1}C_{m-1}^i h^i x_{i+1}=\sum_{i=0}^{m-1}C_{m-1}^i h^i(-1)^{m-i}(C_{m+k-i-1}^{m-i}-2)h^{m-i}r$$

$$=(-1)^m h^m r\left[C_k^m-2\sum_{i=0}^{m-1}(-1)^i C_{m-1}^i\right]$$

$$=(-1)^m h^m r\left[C_k^m-2(1-1)^{m-1}\right]$$

$$=(-1)^m C_k^m h^m r \tag{121}$$

2) Similarly, for the points $a_{-k}(k\ge1)$ and $b_{+k}(k\ge2)$, also obtain

$$y=(-1)^{m-1}C_k^m h^m r \tag{122}$$

Thus both 3) and 4) hold. □

The coordinate of point $a_k(k\ge1)$ is

$$x_i=(-1)^{i-1}C_{m+k-i}^{m-i+1}h^{m-i+1}rs,\ 1\le i\le m \tag{123}$$

The coordinate of point $b_k(k\ge2)$ is

$$x_i=(-1)^{i-1}(C_{m+k-i}^{m-i+1}-2)h^{m-i+1}rs,\ 1\le i\le m \tag{124}$$

From Theorem 4.5, it can be seen that $a_0$ (the origin), $a_k(1\le k\le m-1)$ and $b_k(2\le k\le m-1)$ are all on the $(m-1)$-dimensional hyperplane $N_m(m-1)$. Set $(m-1)$-dimensional hyperplane set $N\triangleq\{N_m(k)\}_{k\ge m-1}$. The items in set $N$ are parallel to each other that cut the m-dimensional state space into the mutually independent regions. The following theorem proves that the $(m-1)$-dimensional hyperplane set $N$ is unique, and it only needs to prove that the $(m-1)$-dimensional hyperplane $N_m(m-1)$ is unique.

**Theorem 4.6.** Set the points $a_k(1\le k\le m-1)$ are on the following $(m-1)$-dimensional hyperplane $N_m^{'}(m-1)$.

$$\sum_{i=0}^{m-1}\alpha_{i+1}x_{i+1}=0 \tag{125}$$

Set $\alpha_1=1$, then $\alpha_{i+1}=C_{m-1}^i h^i(0\le i\le m-1)$ and $\alpha_{i+1}(0\le i\le m-1)$ are unique, i.e., $N_m^{'}(m-1)$ is exactly $N_m(m-1)$.

**Proof.** Substitute (123) into (125) and obtain

$$\sum_{i=0}^{m-1}\alpha_{i+1}(-1)^i C_{m+k-i-1}^{m-i}h^{m-i}rs=0 \tag{126}$$

$$\Rightarrow\sum_{i=0}^{m-1}\alpha_{i+1}(-1)^i C_{m+k-i-1}^{m-i}h^{-i}=0 \tag{127}$$

$$\Rightarrow\sum_{i=1}^{m-1}\alpha_{i+1}(-1)^i C_{m+k-i-1}^{m-i}h^{-i}=-C_{m+k-1}^m \tag{128}$$

Take $k=1,2,\cdots,m-1$ in (128) and (128) is written in matrix form.

$$T_0T_1\alpha=\Gamma \tag{129}$$

Where,

$$T_0=\begin{bmatrix}1&1&1&&1\\C_2^1&C_3^2&C_4^3&\cdots&C_m^{m-1}\\C_3^1&C_4^2&C_5^3&&C_{m+1}^{m-1}\\\vdots&&&\ddots&\vdots\\C_{m-1}^1&C_m^2&C_{m+1}^3&\cdots&C_{2m-3}^{m-1}\end{bmatrix} \tag{130}$$

$$T_1=diag\left[(-1)^{m-1}h^{-(m-1)},(-1)^{m-2}h^{-(m-2)},\cdots,-h^{-1}\right] \tag{131}$$

$$\alpha=[\alpha_m,\alpha_{m-1},\cdots,\alpha_2]^T \tag{132}$$

$$\Gamma=-\left[C_m^m,C_{m+1}^m,\cdots,C_{2m-2}^m\right]^T \tag{133}$$

Since $detT_1\ne0$, we only need to prove that $detT_0\ne0$ and so can prove that $\alpha$ has a unique solution. Combined with the points $a_k(1\le k\le m-1)$ are on the $(m-1)$-dimensional hyperplane $N_m(m-1)$, there is that $\alpha_{i+1}=C_{m-1}^i h^i(0\le i\le m-1)$ and $\alpha_{i+1}(0\le i\le m-1)$ are unique.

Because the matrix $T_0$ performs the row elementary transformation without changing its reversibility, the series of row elementary transformations are performed as follows.

$$T_0=\begin{bmatrix} 1 & 1 & 1 & & 1 \\ C_2^1 & C_3^2 & C_4^3 & \cdots & C_m^{m-1} \\ C_3^1 & C_4^2 & C_5^3 & & C_{m+1}^{m-1} \\ \vdots & & & \ddots & \vdots \\ C_{m-1}^1 & C_m^2 & C_{m+1}^3 & \cdots & C_{2m-3}^{m-1} \end{bmatrix}$$

$$\Rightarrow \begin{bmatrix} 1 & 1 & 1 & & 1 \\ 2 & 3 & 4 & \cdots & m \\ 3\cdot 2 & 4\cdot 3 & 5\cdot 4 & & m(m+1) \\ \vdots & & & \ddots & \vdots \\ C_{m-1}^1 & C_m^2 & C_{m+1}^3 & \cdots & C_{2m-3}^{m-1} \end{bmatrix}$$

$$\Rightarrow \begin{bmatrix} 1 & 1 & 1 & & 1 \\ 2 & 3 & 4 & \cdots & m \\ 2^2 & 3^2 & 4^2 & & m^2 \\ \vdots & & & \ddots & \vdots \\ C_{m-1}^{m-2} & C_m^{m-2} & C_{m+1}^{m-2} & \cdots & C_{2m-3}^{m-2} \end{bmatrix}$$

$$\Rightarrow \begin{bmatrix} 1 & 1 & 1 & & 1 \\ 2 & 3 & 4 & \cdots & m \\ 2^2 & 3^2 & 4^2 & & m^2 \\ \vdots & & & \ddots & \vdots \\ 2^{m-2} & 3^{m-2} & 4^{m-2} & \cdots & m^{m-2} \end{bmatrix}$$

$$\Rightarrow \begin{bmatrix} 1 & 1 & 1 & & 1 \\ 1 & 2 & 3 & \cdots & m-1 \\ 1^2 & 2^2 & 3^2 & & (m-1)^2 \\ \vdots & & & \ddots & \vdots \\ 1^{m-2} & 2^{m-2} & 3^{m-2} & \cdots & (m-1)^{m-2} \end{bmatrix}$$

$$\Rightarrow \begin{bmatrix} 1 & 1 & 1 & & 1 \\ 1 & 2 & 3 & \cdots & m-1 \\ 1 & 2^2 & 3^2 & & (m-1)^2 \\ \vdots & & & \ddots & \vdots \\ 1 & 2^{m-2} & 3^{m-2} & \cdots & (m-1)^{m-2} \end{bmatrix} \triangleq T_0' \tag{134}$$

$$detT_0 \neq 0 \Leftrightarrow detT_0' \neq 0 \tag{135}$$

Then $T_0'$ is exactly the Vandermonde matrix.

$$detT_0' = \prod_{1\le j<i\le m-1}(i-j)\neq 0 \tag{136}$$

$$\Rightarrow detT_0 \neq 0 \tag{137}$$

□

*4.3. (m-1)-dimensional hyperplane $M_m(k)$*

**Definition 4.7.** $k \ge m-1$, the function is defined by

$$z(k) \triangleq \sum_{i=0}^{m-1} C_k^i h^i x_{i+1} \tag{138}$$

**Definition 4.8.** $k \ge m-1$, the $(m-1)$-dimensional hyperplane $M_m(k)$ between $N_m(k)$ and $N_m(k-1)$ is defined by

$$z(k) = (-1)^{m-1} C_k^m h^m rs \tag{139}$$

When $k \ge m$, $z(k) > 0$, $M_m(k)$ is denoted as $M_m(+k)$; when $k \ge m$, $z(k) < 0$, $M_m(k)$ is denoted as $M_m(-k)$; when $k = m-1$, $M_m(k)$ is denoted as $M_m(m-1)$, $z(k)$ equals zero, and then $M_m(m-1)$ is exactly $N_m(m-1)$.

**Theorem 4.9.** $k \ge m-1$, the points $a_k$ and $a_{k-1}$ are on the $(m-1)$-dimensional hyperplane $M_m(k)$.

**Proof.** Substitute the coordinate of $a_k$ (123) into (138) and notice (31) in Theorem 2.16.

$$z(k) = \sum_{i=0}^{m-1} C_k^i h^i (-1)^i C_{m+k-i-1}^{m-i} h^{m-i} rs$$

$$= h^m rs \sum_{i=0}^{m-1} (-1)^i C_k^i C_{m+k-i-1}^{m-i} = (-1)^{m-1} C_k^m h^m rs \tag{140}$$

The coordinate of $a_{k-1}$ is

$$x_i = (-1)^{i-1} C_{m+k-i-1}^{m-i+1} h^{m-i+1} rs, \quad 1 \le i \le m \tag{141}$$

Substitute (141) into (138) and notice (32) in Theorem 2.16.

$$z(k) = \sum_{i=0}^{m-1} C_k^i h^i (-1)^i C_{m+k-i-2}^{m-i} h^{m-i} rs$$

$$= h^m rs \sum_{i=0}^{m-1} (-1)^i C_k^i C_{m+k-i-2}^{m-i} = (-1)^{m-1} C_k^m h^m rs \tag{142}$$

□

**Definition 4.10.** $k \ge m-1$, the $(m-1)$-dimensional hyperplane $\bar{M}_m(k)$ between $N_m(k)$ and $N_m(k-1)$ is defined by

$$z(k) = (-1)^{m-1} (C_k^m - 2C_{k-1}^{m-1}) h^m rs \tag{143}$$

When $k \ge m$, $z(k) > 0$, $\bar{M}_m(k)$ is denoted as $\bar{M}_m(+k)$; when $k \ge m$, $z(k) < 0$, $\bar{M}_m(k)$ is denoted as $\bar{M}_m(-k)$ ; when $k=m-1$ , $\bar{M}_m(k)$ is denoted as $\bar{M}_m(m-1)$ , $z(k)$ equals zero, and then $\bar{M}_m(m-1)$ is exactly $N_m(m-1)$.

**Theorem 4.11.** $k \ge m-1$, the points $b_k$ and $b_{k-1}$ are on the $(m-1)$ -dimensional hyperplane $\bar{M}_m(k)$.

**Proof.** Substitute the coordinate of $b_k$ (124) into (138) and notice Theorem 2.17 and (31) in Theorem 2.16.

$$z(k)=\sum_{i=0}^{m-1} C_k^i h^i (-1)^i (C_{m+k-i-1}^{m-i}-2)h^{m-i}rs$$

$$=h^m rs\left[(-1)^{m-1} C_k^m - 2\sum_{i=0}^{m-1}(-1)^i C_k^i\right]$$

$$=(-1)^{m-1}(C_k^m - 2C_{k-1}^{m-1})h^m rs \tag{144}$$

The coordinate of $b_{k-1}$ is

$$x_i=(-1)^{i-1}(C_{m+k-i-1}^{m-i+1}-2)h^{m-i+1}rs,\ \ 1\le i\le m \tag{145}$$

Substitute (145) into (138) and notice Theorem 2.17 and (32) in Theorem 2.16.

$$z(k)=\sum_{i=0}^{m-1} C_k^i h^i (-1)^i (C_{m+k-i-2}^{m-i}-2)h^{m-i}rs$$

$$=h^m rs\left[(-1)^{m-1} C_k^m - 2\sum_{i=0}^{m-1}(-1)^i C_k^i\right]$$

$$=(-1)^{m-1}(C_k^m - 2C_{k-1}^{m-1})h^m rs \tag{146}$$

□

Obviously, the $(m-1)$ -dimensional hyperplanes $\bar{M}_m(k)$ and $M_m(k)$ are parallel to each other. This means that $\bar{M}_m(k)$ between $N_m(k)$ and $N_m(k-1)$ is the $(m-1)$-dimensional hyperplane passing through the points $b_k, b_{k-1}$ and also parallel to $M_m(k)$.

**Definition 4.12.** $k \ge m-1$ , the $(m-1)$ -dimensional hyperplane $M_m^\beta(k)$ between $N_m(k)$ and $N_m(k-1)$ is defined by

$$z(k)=(-1)^{m-1}[C_k^m - 2(1-\beta)C_{k-1}^{m-1}]h^m rs \tag{147}$$

When $k\ge m$ , $z(k)>0$ , $M_m^\beta(k)$ is denoted as $M_m^\beta(+k)$; when $k\ge m$, $z(k)<0$, $M_m^\beta(k)$ is denoted as $M_m^\beta(-k)$ ; when $k=m-1$ , $M_m^\beta(k)$ is denoted as $M_m^\beta(m-1)$ , $z(k)$ equals zero, and then $M_m^\beta(m-1)$ is exactly $N_m(m-1)$.

The following Theorem 4.13 can be easily obtained from Theorem 4.9 and Theorem 4.11.

**Theorem 4.13.** $k \ge m-1$ , $0\le\beta\le 1$ the point $X_\beta = \beta a_k + (1-\beta)b_k$ on the line segment $a_k b_k$ $(a_{+k}b_{-k}\ or\ a_{-k}b_{+k})$ is on the $(m-1)$ -dimensional hyperplane $M_m^\beta(k)$.

Obviously, when $\beta=1$ , the point $X_\beta$ is $a_k$ and $M_m^\beta(k)$ is exactly $M_m(k)$; when $\beta=0$, the point $X_\beta$ is $b_k$ and $M_m^\beta(k)$ is precisely $\bar{M}_m(k)$ . The $(m-1)$ -dimensional hyperplanes $M_m^\beta(k)$ and $M_m(k)$ are parallel to each other, and this means that $M_m^\beta(k)$ between $N_m(k)$ and $N_m(k-1)$ is the $(m-1)$-dimensional hyperplane passing through the point $X_\beta$ and parallel to $M_m(k)$.

## *4.4. m-order synthesis function of nonlinear region*

Suppose that $k\ge m$ , $m\ge 2$ , the state $X \triangleq [x_1,x_2,\cdots,x_m]^T$ is located between the parallel hyperplanes $N_m(k)$ and $N_m(k-1)$, and satisfies the following equation.

$$y=\sum_{i=0}^{m-1} C_{m-1}^i h^i x_{i+1} = C_{k'}^m h^m rs \tag{148}$$

Where, $k-1<k'\le k$ , $k$ is a positive integer and meets the following equations from (148).

$$\begin{cases}\prod_{i=0}^{m-1}(k'-i)=\dfrac{m!|y|}{h^m r}\\ k-1<k'\le k\\ k\ge m\end{cases} \tag{149}$$

The value $k$ can be obtained from (149).

Now the control target has 2 points: the state $X$ is controlled to reach the region between $N_m(k-1)$ and $N_m(k-2)$ from the region between $N_m(k)$ and $N_m(k-1)$; the state $X$ is controlled to reach the state on $M_m(k-1)$ between $N_m(k-1)$ and $N_m(k-2)$. It is assumed that the state on $M_m(k-1)$ between $N_m(k-1)$ and $N_m(k-2)$ is $X' \triangleq [x_1', x_2', \cdots, x_m']^T$ after the state $X$ is controlled by one step. And $M_m(k-1)$ between $N_m(k-1)$ and $N_m(k-2)$ meets the following equation.

$$z(k-1)=\sum_{i=0}^{m-1} C_{k-1}^i h^i x_{i+1}' = (-1)^{m-1} C_{k-1}^m h^m rs \tag{150}$$

The new state $X'$ meets

$$\begin{cases} x_1' = x_1 + hx_2 \\ x_2' = x_2 + hx_3 \\ \vdots \\ x_m' = x_m + hu \end{cases} \tag{151}$$

Substitute (151) into (150), notice (5) and (7) and obtain

$$\sum_{i=0}^{m-1} C_k^i h^i x_{i+1} + C_{k-1}^{m-1} h^m u = (-1)^{m-1} C_{k-1}^m h^m rs \tag{152}$$

$$\Rightarrow u = (-1)^m \left(1 - \frac{k}{m}\right) rs - \frac{\sum_{i=0}^{m-1} C_k^i h^i x_{i+1}}{C_{k-1}^{m-1} h^m} \tag{153}$$

**Definition 4.14.** $k \geq m$, $m \geq 2$, $|y| > h^m r$, the function is defined as follows.

$$a(x_1, x_2, \cdots, x_m, r, h) \triangleq$$

$$(-1)^{m-1} \left(1 - \frac{k}{m}\right) rs + \frac{\sum_{i=0}^{m-1} C_k^i h^i x_{i+1}}{C_{k-1}^{m-1} h^m} \tag{154}$$

Since $|u| \leq r$ and $|y| > h^2 r$, then the time optimal control input (i.e., the m-order synthesis function of nonlinear region) is

$$u = -r \cdot sat\left[a(x_1, x_2, \cdots, x_m, r, h), r\right] \tag{155}$$

Further, if the state $X$ is located on $M_m(k)$ between $N_m(k)$ and $N_m(k-1)$, and the state $X$ satisfies (139). We substitute (139) into (152), notice (7) and obtain

$$(-1)^{m-1} C_k^m h^m rs + C_{k-1}^{m-1} h^m u = (-1)^{m-1} C_{k-1}^m h^m rs \tag{156}$$

$$\Rightarrow u = (-1)^m rs \tag{157}$$

(157) means that the state $X$ can arrive at $M_m(k-1)$ (but not sure is the region between $N_m(k-1)$ and $N_m(k-2)$) by one step control under the control input of $u = (-1)^m rs$.

It can be seen that once the state $X$ is located on $M_m(k)$ between $N_m(k)$ and $N_m(k-1)$, it will reach the someone state on $M_m(k-1)$ under the control of $u = (-1)^m rs$ until it reaches the state on $M_m(m)$, and then enters the linear region (the region of $|y| \leq h^m r$), i.e., reaches the state on $M_m(m-1)$ (or $N_m(m-1)$). In the linear region, the linear control is used to make the state arrive at the origin.

Set $(m-1)$-dimensional hyperplane set $M \triangleq \{M_k\}_{k \geq m-1}$. It is demonstrated below that the hyperplane set $M$ that satisfies the above control requirements is unique. Since Theorem 4.6 guarantees the uniqueness of $N_m(m-1)$ (or $M_m(m-1)$), we prove in turn that all terms in set $M$ are unique by using this as a starting point.

**Theorem 4.15.** $k \geq m-1$, $m \geq 2$, suppose that the $(m-1)$-dimensional hyperplane $\hat{M}_m(k)$ between $N_m(k)$ and $N_m(k-1)$ satisfies the following equation.

$$\hat{z}(k) \overset{def}{=} \sum_{i=0}^{m-1} \alpha_i(k) h^i x_{i+1} = (-1)^{m-1} \beta(k) h^m rs \tag{158}$$

Set $(m-1)$-dimensional hyperplane set $\hat{M} \triangleq \{\hat{M}_k\}_{k \geq m-1}$.

If the following conditions are satisfied:

1) $\hat{z}(m-1) = z(m-1) = 0$ (159)

i.e., $\alpha_i(m-1) = C_{m-1}^i$, $\beta(m-1) = 0$ (160)

So $\hat{M}_m(m-1)$ is exactly $M_m(m-1)$.

2) Set $\alpha_{m-1}(k) = C_k^{m-1}$ (161)

3) The state $X$ on $\hat{M}_m(k+1)$ between $N_m(k+1)$ and $N_m(k)$ reaches the state $X'$ on $\hat{M}_m(k)$ by one step control under the control of $u = (-1)^m rs$, and the state $X'$ satisfies (151).

Then, the hyperplane set $\hat{M}$ is unique, that is, $\hat{M}$ is precisely $M$.

**Proof.** Let $k = m-1$ in the condition 3), and obtain

$$\sum_{i=0}^{m-1} \alpha_i(m-1) h^i x_{i+1}' = 0 \tag{162}$$

Substitute (151) into (162), utilize (160) and notice (5) and (7).

$$\Rightarrow \sum_{i=0}^{m-1} C_m^i h^i x_{i+1} + C_{m-1}^{m-1} h^m u = 0 \tag{163}$$

$$\overset{(157)}{\Rightarrow} \sum_{i=0}^{m-1} C_m^i h^i x_{i+1} = (-1)^{m-1} C_{m-1}^{m-1} h^m rs \tag{164}$$

Compare (158) and (164), notice (161) and obtain

$$\alpha_i(m) = C_m^i,\ \beta(m) = C_{m-1}^{m-1} = C_m^m \tag{165}$$

Thereby, $\hat{M}_m(m)$ is exactly $M_m(m)$.

Then set $k = m$ in the condition 3), and obtain

$$\sum_{i=0}^{m-1} \alpha_i(m) h^i x_{i+1}' = (-1)^{m-1} \beta(m) h^m rs \tag{166}$$

Substitute (151) into (166), utilize (165) and notice (5) and (7).

$$\Rightarrow \sum_{i=0}^{m-1} C_{m+1}^{i} h^{i} x_{i+1} + C_{m}^{m-1} h^{m} u = (-1)^{(m-1)} C_{m}^{m} h^{m} rs \quad (167)$$

$$\overset{(157)}{\Rightarrow} \sum_{i=0}^{m-1} C_{m+1}^{i} h^{i} x_{i+1} = (-1)^{m-1} \left( C_{m}^{m-1} + C_{m}^{m} \right) h^{m} rs$$

$$= (-1)^{m-1} C_{m+1}^{m} h^{m} rs \quad (168)$$

Compare (158) and (168), notice (161) and obtain

$$\alpha_i(m+1) = C_{m+1}^{i}, \ \beta(m+1) = C_{m+1}^{m} \quad (169)$$

This means that $\hat{M}_m(m+1)$ is precisely $M_m(m+1)$.

And so on, we get by recursion in turn that $\hat{M}_m(m+2)$ is $M_m(m+2)$, $\hat{M}_m(m+3)$ is $M_m(m+3), \cdots$. Thus, for $k \ge m-1$, $\hat{M}_m(k)$ is $M_m(k)$, and there are

$$\alpha_i(k) = C_k^{i}, \ \beta(k) = C_k^{m} \quad (170)$$

So the hyperplane set $\hat{M}$ is unique, that is, $\hat{M}$ is precisely $M$. □

The above analysis shows that the piecewise combination of all $M_m(k)(k \ge m)$ regions between $N_m(k)$ and $N_m(k-1)$ is a kind of $(m-1)$-dimensional time optimal trajectory to reach the origin in nonlinear region of m-order discrete system. The conclusion is consistent with the case of $m = 2$.

**Theorem 4.16.** $0 \le \beta \le 1$, $m \ge 2$, $k \ge m$, suppose that the state $X$ is located on $M_m^{\beta}(k)$ between $N_m(k)$ and $N_m(k-1)$, $X$ is controlled by one step to reach the state $X'$ on $M_m(k-1)$ (but not sure is the region between $N_m(k-1)$ and $N_m(k-2)$), and then the control input of state $X$ is

$$u = (-1)^{m} (2\beta - 1) rs \quad (171)$$

**Proof.** Substitute (151) into (150), utilize (5) and (7) to obtain (152), and then substitute (147) into (152).

$$(-1)^{m-1} \left[ C_k^{m} - 2(1-\beta) C_{k-1}^{m-1} \right] h^{m} rs + C_{k-1}^{m-1} h^{m} u$$

$$= (-1)^{m-1} C_{k-1}^{m} h^{m} rs \quad (172)$$

$$\Rightarrow C_{k-1}^{m-1} u = (-1)^{m} (2\beta - 1) C_{k-1}^{m-1} rs \quad (173)$$

□

In Theorem 4.16, the state $X$ is located on $M_m(k)$ between $N_m(k)$ and $N_m(k-1)$, and $u = (-1)^{m} rs$ if $\beta = 1$; $X$ is located on $\bar{M}_m(k)$ between $N_m(k)$ and $N_m(k-1)$, and $u = (-1)^{m-1} rs$ if $\beta = 0$; $X$ is located on $M_m^{0.5}(k)$ between $N_m(k)$ and $N_m(k-1)$, and $u = 0$ if $\beta = 0.5$. This means that in the region between $N_m(k)$ and $N_m(k-1)$, the control input takes the value of $(-1)^{m}(2\beta-1)rs$ between $(-1)^{m} rs$ and $(-1)^{m-1} rs$ if the state $X$ is located in the region section between $M_m(k)$ and $\bar{M}_m(k)$; the control takes $(-1)^{m} rs$ in the region section outside $M_m(k)$; and then the control takes $(-1)^{m-1} rs$ in the region section outside $\bar{M}_m(k)$. This feature is consistent with the case of $m = 2$.

**Theorem 4.17.** $0 \le \beta \le 1$, $m \ge 2$, $k \ge m$, suppose that the state $X$ is located on the line segment $a_k a_{k-1}$, $X = \beta a_k + (1-\beta) a_{k-1}$, $X$ is controlled by one step to reach the state $X'$ under the control of $u = (-1)^{m} rs$, and then it is easy to obtain that $X' = \beta a_{k-1} + (1-\beta) a_{k-2}$ from the linear properties of discrete system (89), i.e., the state $X'$ is located on the line segment $a_{k-1} a_{k-2}$.

**Definition 4.18.** The $(m-1)$-dimensional hyperplane $\bar{N}_m(k)$ is defined to meet the following equation.

$$\bar{y} = \bar{y}(X) \triangleq \sum_{i=0}^{m-1} C_m^{i} h^{i} x_{i+1} = \left[ C_{k-1}^{m} + (-1)^{m-1} \right] h^{m} rs \quad (174)$$

When $k = m, m-1, \cdots, 1$, $\bar{N}_m(k)$ is the same hyperplane uniformly denoted as $\bar{N}_m(m)$.

**Theorem 4.19.** $a_k (k \ge 1)$ is located on $\bar{N}_m(k)$, but $b_k (k \ge 2)$ is not located on $\bar{N}_m(k)$.

**Proof.** For $a_k (k \ge 1)$, substitute (123) into (174), utilize Theorem 2.18 and obtain

$$\bar{y} \triangleq \sum_{i=0}^{m-1} C_m^{i} h^{i} x_{i+1} = \sum_{i=0}^{m-1} C_m^{i} h^{i} (-1)^{i} C_{m+k-i-1}^{m-i} h^{m-i} rs$$

$$= h^{m} rs \sum_{i=0}^{m-1} (-1)^{i} C_m^{i} C_{m+k-i-1}^{m-i} = \left[ C_{k-1}^{m} + (-1)^{m-1} \right] h^{m} rs \quad (175)$$

For $b_k (k \ge 2)$, substitute (124) into (174), utilize Theorem 2.18 and Theorem 2.19, and obtain

$$\bar{y} \triangleq \sum_{i=0}^{m-1} C_m^{i} h^{i} x_{i+1} = \sum_{i=0}^{m-1} C_m^{i} h^{i} (-1)^{i} \left( C_{m+k-i-1}^{m-i} - 2 \right) h^{m-i} rs$$

$$= \left[ C_{k-1}^{m} + (-1)^{m-1} - 2 \sum_{i=0}^{m-1} (-1)^{i} C_m^{i} \right] h^{m} rs$$

$$= \left[ C_{k-1}^{m} + (-1)^{m-1} - 2(-1)^{m-1} \right] h^{m} rs$$

$$=\left[C_{k-1}^m+(-1)^m\right]h^m rs \tag{176}$$

And so, $a_k\ (k\ge 1)$ is located on $\bar{N}_m(k)$, but $b_k\ (k\ge 2)$ is not located on $\bar{N}_m(k)$. □

**Theorem 4.20.** $k\ge m+1$, suppose that the state $X$ is located on $M_m(k)$ between $N_m(k)$ and $N_m(k-1)$, and $X$ is controlled by one step to reach the state $X'$ on $M_m(k-1)$ between $N_m(k-1)$ and $N_m(k-2)$ under the control of $u=(-1)^m rs$. The point $\bar{X}$ is the intersection point of the line segment $a_k a_{k-1}$ and hyperplane $\bar{N}_m^X(k)$ passing through the point $X$ and parallel to $\bar{N}_m(k)$, and the point $\bar{X}'$ is the intersection point of the line segment $a_{k-1}a_{k-2}$ and hyperplane $N_m^{X'}(k-1)$ passing through the point $X'$ and parallel to $N_m(k-1)$. $0<\beta'\le 1$, Set

$$\bar{X}=\beta a_k+(1-\beta)a_{k-1} \tag{177}$$

$$\bar{X}'=\beta' a_{k-1}+(1-\beta')a_{k-2} \tag{178}$$

Then $\beta=\beta'$.

**Proof.**

$$y(X')=y(\bar{X}')=\beta' y(a_{k-1})+(1-\beta')y(a_{k-2})$$

$$=\left[\beta' C_{k-1}^m+(1-\beta')C_{k-2}^m\right]h^m rs \tag{179}$$

$$\bar{y}(X)=\bar{y}(\bar{X})$$

$$=\bar{y}\left(\beta a_k+(1-\beta)a_{k-1}\right)=\beta\bar{y}(a_k)+(1-\beta)\bar{y}(a_{k-1})$$

$$=\left[\beta\left(C_{k-1}^m+(-1)^{m-1}\right)+(1-\beta)\left(C_{k-2}^m+(-1)^{m-1}\right)\right]h^m rs$$

$$=\left[\beta C_{k-1}^m+(1-\beta)C_{k-2}^m+(-1)^{m-1}\right]h^m rs \tag{180}$$

Utilize (151), (157) and (180), notice (5) and (7),

$$y(X')=\sum_{i=0}^{m-1}C_{m-1}^i h^i x_{i+1}'$$

$$=\sum_{i=0}^{m-1}C_{m-1}^i h^i x_{i+1}+\sum_{i=0}^{m-2}C_{m-1}^i h^{i+1}x_{i+2}+C_{m-1}^{m-1}h^m u$$

$$=\sum_{i=0}^{m-1}C_m^i h^i x_{i+1}+(-1)^m h^m rs=\bar{y}(X)+(-1)^m h^m rs$$

$$=\left[\beta C_{k-1}^m+(1-\beta)C_{k-2}^m\right]h^m rs \tag{181}$$

Obtain from (179) and (181)

$$(\beta-\beta')(C_{k-1}^m-C_{k-2}^m)=0 \tag{182}$$

For $k\ge m+1$, $C_{k-1}^m-C_{k-2}^m\ne 0$ (183)

$$\Rightarrow \beta-\beta'=0 \tag{184}$$

□

The region of $M_m(k)$ between $N_m(k)$ and $N_m(k-1)$ is denoted as $\Phi_1(k)$, the region of $M_m(k)$ between $\bar{N}_m(k)$ and $\bar{N}_m(k-1)$ is denoted as $\Phi_2(k)$, and $\Phi(k)\triangleq\Phi_1(k)\cap\Phi_2(k)$. Then $\Phi(k)\subset\Phi_1(k)$, $k\ge m+1$, $m\ge 3$; when $k=m$, $m\ge 3$, $\bar{N}_m(m)$ is the same hyperplane as $\bar{N}_m(m-1)$, but $N_m(m)$ is not same as $\bar{N}_m(m)$, and so $\Phi(k)\subset\Phi_1(k)$. Thence, for $k\ge m$, $m\ge 3$, $\Phi(k)\subset\Phi_1(k)$.

Theorem 4.20 means that the region that the state $X$ is located is $\Phi(k)$ if the state $X'$ is located on $M_m(k-1)$ between $N_m(k-1)$ and $N_m(k-2)$, i.e., some states in the region $\Phi_1(k)$ that do not belong to the region $\Phi(k)$ need at least two steps to reach the region $\Phi_1(k-1)$ under the control of $u=(-1)^m rs$. This property is established only when $k\ge m$ and $m\ge 3$, which is not the same as the case of $m=2$.

### *4.5. m-order synthesis function of linear region*

In order to derive the m-order synthesis function in linear region of $M_m(m-1)$, a series of hyperplanes are further defined in $M_m(m-1)$.

**Definition 4.21.** $m\ge 2$, $0\le\nu\le m-1$, $0\le k\le m-1-\nu$, the function is defined as follows.

$$z_{m-\nu}(m-\nu-1)\triangleq\sum_{i=0}^{m-(\nu+1)}C_{m-(\nu+1)}^i h^i x_{i+(\nu+1)} \tag{185}$$

In $M_m(m-1)$, the $(m-1-\nu)$-dimensional hyperplane $M_{m-\nu}(m-1-\nu)$ is defined to satisfy the following equations.

$$z_{m-\mu}(m-\mu-1)=0,\ \ \mu=0,1,\cdots,\nu \tag{186}$$

Thus, Definition 4.21 defines a series of hyperplanes in $M_m(m-1)$: $M_{m-1}(m-2)$, $M_{m-2}(m-3)$, $\cdots$, $M_2(1)$, $M_1(0)$. And $M_1(0)$ is the origin. Then there is the following theorem.

**Theorem 4.22.** If $m\ge 2$, $0\le\nu\le m-1$, $0\le k\le m-1-\nu$, then

1) The points $\{a_k\}_{0\le k\le m-1-\nu}$ and $\{b_k\}_{2\le k\le m-1-\nu}$ are all located on $M_{m-\nu}(m-1-\nu)$;

2) Obviously, $M_m(m-1) \supset M_{m-1}(m-2) \supset \cdots \supset M_2(1) \supset M_1(0)$;

3) It is as easy to prove as Theorem 4.6: the set $\{M_{m-\nu}(m-1-\nu)\}_{0\le\nu\le m-1}$ is unique.

**Proof.** We only prove 1) as follows.

Obviously, the point $a_0$ satisfies (186).

Substitute the coordinate of $a_k\,(1\le k\le m-1-\nu)$ (123) into (185), use Theorem 2.20 and obtain

$\mu = 0,1,\cdots,\nu$

$$\sum_{i=0}^{m-(\mu+1)} C^i_{m-(\mu+1)} h^i x_{i+(\mu+1)}$$

$$= \sum_{i=0}^{m-(\mu+1)} C^i_{m-(\mu+1)} h^i (-1)^{i+\mu} C^{m-(\mu+i)}_{k+m-(\mu+i)-1} h^{m-(\mu+i)} rs$$

$$= (-1)^{\mu} h^{m-\mu} rs \sum_{i=0}^{m-(\mu+1)} (-1)^i C^i_{m-(\mu+1)} C^{m-(\mu+i)}_{k+m-(\mu+i)-1} = 0 \quad (187)$$

$\Rightarrow$ The point $a_k\,(1\le k\le m-1-\nu)$ is all located on $M_{m-\nu}(m-1-\nu)$.

Substitute the coordinate of $b_k\,(2\le k\le m-1-\nu)$ (124) into (185), utilize Theorem 2.20 and obtain

$\mu = 0,1,\cdots,\nu$

$$\sum_{i=0}^{m-(\mu+1)} C^i_{m-(\mu+1)} h^i x_{i+(\mu+1)}$$

$$= (-1)^{\mu} h^{m-\mu} rs \sum_{i=0}^{m-(\mu+1)} (-1)^i C^i_{m-(\mu+1)} [C^{m-(\mu+i)}_{k+m-(\mu+i)-1} - 2]$$

$$= 2(-1)^{\mu} h^{m-\mu} rs \sum_{i=0}^{m-(\mu+1)} (-1)^i C^i_{m-(\mu+1)}$$

$$= 2(-1)^{\mu} h^{m-\mu} rs (1-1)^{m-(\mu+1)} = 0 \quad (188)$$

$\Rightarrow$ The point $b_k\,(2\le k\le m-1-\nu)$ is all located on $M_{m-\nu}(m-1-\nu)$. □

The state $X$ between $N_m(m)$ and $N_m(m-1)$ is controlled to reach the state $X'$ on $M_m(m-1)$ which satisfies the following equation.

$$\sum_{i=0}^{m-1} C^i_{m-1} h^i x'_{i+1} = 0 \quad (189)$$

Substitute (151) into (189), utilize (5) and (7), and obtain

$$\sum_{i=0}^{m-1} C^i_m h^i x_{i+1} + C^{m-1}_{m-1} h^m u = 0 \quad (190)$$

$$\Rightarrow u = -\frac{\sum_{i=0}^{m-1} C^i_m h^i x_{i+1}}{h^m} \quad (191)$$

Further, for $m\ge 2$, $\nu = 0,1,\cdots,m-1$, the state $X$ on $M_{m-\nu}(m-1-\nu)$ is controlled by one step to reach the state $X'$ on $M_{m-\nu-1}(m-2-\nu)$ which satisfies the following equations.

$$\sum_{i=0}^{m-(\mu+1)} C^i_{m-(\mu+1)} h^i x'_{i+(\mu+1)} = 0,\ \ \mu = 0,1,\cdots,\nu \quad (192)$$

Substitute (151) into (192), use (5) and (7), and then obtain

$$\sum_{i=0}^{m-(\mu+1)} C^i_{m-\mu} h^i x_{i+(\mu+1)} + C^{m-(\mu+1)}_{m-(\mu+1)} h^{m-\mu} u = 0, \mu = 0,1,\cdots,\nu \quad (193)$$

$$\Rightarrow u = -\frac{\sum_{i=0}^{m-(\mu+1)} C^i_{m-\mu} h^i x_{i+(\mu+1)}}{h^{m-\mu}},\ \ \mu = 0,1,\cdots,\nu \quad (194)$$

For the state $X$ on $M_{m-\nu}(m-1-\nu)$, $\mu = 0$ is selected uniformly, that is, the uniform control input expression (191) is used in the following.

Thus, in the linear region of $|y|\le h^m r$, the state $X$ on $M_m(m-1)$ is controlled under the expression of (191) to gradually reach $M_{m-1}(m-2)$, $M_{m-2}(m-3)$, $\cdots$, $M_2(1)$, $M_1(0)$, that is, finally arrive at the origin.

The nested combination of $M_m(m-1)$, $M_{m-1}(m-2)$, $\cdots$, $M_2(1)$, $M_1(0)$, is a kind of time optimal trajectory to reach the origin in linear region of m-order discrete system.

For $|u|\le r$ and $|y|\le h^m r$, the time optimal control input (i.e., the m-order synthesis function of linear region) is

$$u = -r\cdot sat\left(\frac{\sum_{i=0}^{m-1} C^i_m h^i x_{i+1}}{h^m}, r\right) \quad (195)$$

### *4.6. m-order synthesis function*

**Definition 4.23.** In both cases of $|y| > h^m r$ and $|y|\le h^m r$, the function $a(x_1,x_2,\cdots,x_m,r,h)$ is uniformly defined as follows.

$$a(x_1,x_2,\cdots,x_m,r,h) \triangleq \begin{cases} \dfrac{\sum_{i=0}^{m-1} C^i_m h^i x_{i+1}}{h^m}, & |y|\le h^m r \\ (-1)^{m-1}\left(1-\dfrac{k}{m}\right) rs + \dfrac{\sum_{i=0}^{m-1} C^i_k h^i x_{i+1}}{C^{m-1}_{k-1} h^m}, & |y| > h^m r \end{cases} \quad (196)$$

Where, $k$ is a positive integer and satisfies (149).

**Theorem 4.24.** The time optimal control input of m-order discrete system (89) (denoted as $fxiao(x_1, x_2, \cdots, x_m, r, h)$) is

$$u = fxiao(x_1, x_2, \cdots, x_m, r, h) \triangleq -r \cdot sat\left[a(x_1, x_2, \cdots, x_m, r, h), r\right] \tag{197}$$

(197) is the m-order synthesis function of discrete system (89) with the origin as the end point. Obviously, this is a nonlinear function.

**Theorem 4.25.** If the control input limit $|u| \le r$ is removed in the m-order discrete system (89), i.e., $r \to +\infty$, then we can obtain the m-order linear synthesis function with the origin as the end point as follows.

$$u = -\frac{\sum_{i=0}^{m-1} C_m^i h^i x_{i+1}}{h^m} \tag{198}$$

*4.7. m-order synthesis function of tracking problem*

Any real feedback control is achieved by detecting feedback errors, and "eliminating errors based on feedback errors" is the most basic idea and method of feedback control. Any feedback control also has a control target, which can be a given constant or a time-varying trajectory, and the former refers to the adjustment problem, while the latter is called the tracking problem. From the error equation, the adjustment and tracking problems are only the distinction between the constant disturbance and time-varying disturbance. However, it is not necessary to distinguish them from the disturbance suppression effect and they can be treated uniformly and indiscriminately as the tracking problem (Han, 2008).

The m-order synthesis function of discrete system (89) based on the tracking problem is given in the following, and for the theoretical basis on which it is based, please refer to the general theoretical conclusions and their proofs in (Han, 2008; Han & Wang, 1994).

**Theorem 4.26.** For the given signal $v(t)$ (abbreviated as $v$), which can be a given constant or a time-varying trajectory, the m-order synthesis function of discrete system (89) with $v$ as the control target is

$$u = fxiao(x_1 - v, x_2, \cdots, x_m, r, h) \tag{199}$$

(199) ensures that the approximated differential signals of $v$ from 1st-order to $(m-1)$-order are extracted, while the system (89) tracks the given signal as quickly as possible.

**Definition 4.27.** For the integer $\mu \ge 0$, $\mu$-order differential operator is denoted as $D^\mu \triangleq \frac{d^\mu}{dt^\mu}$. The operator $D^0$ is the original signal, and for the given signal $v$, $D^0 v = v$.

Then, for $0 \le \mu \le m-1$, $x_{\mu+1}$ is the approximate signal of $D^\mu v$, and the smaller the control step length $h$ is, the more the $x_{\mu+1}$ approaches the signal $D^\mu v$. However, there is also a clear flaw: the time optimal control system do not has good filtering performance when the system is affected by random disturbances. Because the response to input is too fast, the input noise is reproduced as much as possible in the output signal (or feedback signal). This shows that there is a contradiction between the rapid demand for time optimal control and good noise immunity requirement for control system (Qian & Song, 2011). Han Jingqing has pioneered the research on how to take a compromise between the two to gain comprehensive control effects (Han, 2008). The basic idea of Han Jingqing is used to solve this problem in the following.

## 5. Noise suppression and signal extraction

For system (89), the variable $h$ in m-order synthesis function (199) is changed to a new variable $h_0$ that is independent of control step $h$, and take $h_0 \triangleq n_0 h$, $n_0 > 1$. By increasing the value of $n_0$ to provide sufficiently strong filtering performance, the effect of noise can be effectively suppressed (Han, 2008).

**Theorem 5.1.** For the given signal $v$, the m-order tracking-form synthesis function of discrete system (89) with $v$ as the control target is

$$u = fxiao(x_1 - v, x_2, \cdots, x_m, r, n_0 h) \tag{200}$$

(200) not only guarantees the time optimal control requirement of system (89), but also takes into account the noise suppression requirement of control system. The case of m-order linear synthesis function is only analyzed to illustrate the key role of value $n_0$ in the following.

In (200), the smaller the control step $h$ is, the more the discrete signal $x_{\mu+1}(k)(0 \le \mu \le m-1)$ approaches the continuous signal $D^\mu v$, and in this case, the discrete system can be restored to the m-order continuous system to be analyzed. The linear time optimal control continuous system with $v$ as the control target is

$$\begin{cases} u = -\frac{(x_1 - v) + \sum_{i=1}^{m-1} C_m^i (n_0 h)^i (D^i x_1)}{(n_0 h)^m} \\ \dot{x}_1 = x_2 \\ \dot{x}_2 = x_3 \\ \vdots \\ \dot{x}_m = u \end{cases} \tag{201}$$

Considering the zero-state response of system, the Laplace operator $s$ is used instead of the differential operator $D$, and the transfer relation is obtained as follows from (201).

$$x_{i+1} = \frac{s^i}{(1 + n_0 hs)^m} v, \quad 0 \le i \le m-1 \tag{202}$$

In (202), the relative order of transfer function of $x_{i+1}$ with respect to $v$ is $m-i$, and since $m-i \geq 1$, every $x_{i+1}$ is physically achievable. At the same time, every $x_{i+1}$ has the ability to suppress noise, and further, the larger the value of $m-i$ is, the stronger its ability to suppress noise.

In (202), the m-order filter (composed of m identical first-order filter links in series) has the equivalent m-order integral effect and plays the role of noise suppression. $n_0 h$ is the filter time constant of m-order filter link. The variable $n_0 h$ is called the filter factor in (Han, 2008). But in this paper, $n_0$ is proposed as the filter factor in the following, which is more in line with its physical meaning. While $x_{i+1}$ is the extracted signal obtained after $D^i v$ is filtered by m-order filter. Generally, choosing the value $n_0$ greater than one can play a better filtering effect.

The above analysis shows that the m-order tracking-form synthesis function with filter factor can be regarded as the expansion or development of classical PID, which has the following characteristics:

1) The control strategy of eliminating errors based on feedback errors is derived from the classical PID.

2) The closed-loop control is realized by the state feedback control law based on the tracking-form, which comes from modern control theory.

3) The filter factor is used to effectively suppress the noise's influence, and obtain a better filtering effect (with equivalent integral effect), which stems from modern signal processing technology.

By use of the filter factor, the better filtering effect can be obtained. But at the same time, it also brings the shortcoming of phase delay (or time delay) and amplitude attenuation for the extracted signal and its differential signals. The predictive compensation problem of extracted signals is handled in the following.

## 6. Predictive compensation

The analysis of linear time optimal control closed-loop system with filter factor shows that the extracted signal and its differential signals all have a phase delay of about $mn_0 h$ (when $h$ is small enough) and a certain degree of amplitude attenuation (related to the frequency characteristic of $v$). Since the extracted signal and its differential signals have been filtered, they can be used to predict and compensate the phase delay and amplitude attenuation so as to obtain the compensated signal and its differential signals (Han, 2008).

For the m-order closed-loop system consisting of system (89) and m-order synthesis function (200), the signals $x_{i+1}\left(0 \leq i \leq m-2\right)$ are predicted and compensated as follows.

$i = 0,1,\cdots,m-2$

$$\hat{x}_{i+1} = \left(1 + \frac{mn_0 h}{m-1-i} D\right)^{m-1-i} x_{i+1} \tag{203}$$

$$= \sum_{\mu=0}^{m-1-i} C_{m-1-i}^{\mu} \left(\frac{mn_0 h}{m-1-i}\right)^{\mu} x_{i+1+\mu} \tag{204}$$

Where $D$ is the differential operator.

$\hat{x}_{i+1}\left(0 \leq i \leq m-2\right)$ is the compensated signal for $x_{i+1}$.

The m-order linear system (201) is only analyzed in the following, considering the zero-state response of (203), the Laplace operator $s$ is used instead of the differential operator $D$, and the transfer relation is obtained from (203).

$$\hat{x}_{i+1} = \left(1 + \frac{mn_0 h}{m-1-i} s\right)^{m-1-i} x_{i+1} \tag{205}$$

Obtain from (202) and (205)

$$\hat{x}_{i+1} = \left(1 + \frac{mn_0 h}{m-1-i} s\right)^{m-1-i} \frac{s^i}{\left(1+n_0 hs\right)^m} v \tag{206}$$

All of the transfer functions of $\hat{x}_{i+1}\left(0 \leq i \leq m-2\right)$ with respect to $v$ in (206) have a relative order of one, so they are physically achievable and have the same ability to suppress noise.

When $h$ is small enough, the 1st-order approximation with respect to $h$ can be done on the right of (205) and (206) as follows.

$$\hat{x}_{i+1} \approx \left[1 + \left(m-1-i\right) \frac{mn_0 h}{m-1-i} s\right] x_{i+1} = \left(1 + mn_0 hs\right) x_{i+1} \tag{207}$$

$$\hat{x}_{i+1} \approx \left(1 + mn_0 hs\right) \frac{s^i}{\left(1 + mn_0 hs\right)} v = s^i v \tag{208}$$

$0 \leq i \leq m\text{-}2$ , $\hat{x}_{i+1}$ is ahead of $x_{i+1}$ by phase of about $mn_0 h$ in (207), and (208) shows that $\hat{x}_{i+1}$ is an approximation of $D^i v$. This shows that the phase delay and amplitude attenuation of $x_{i+1}$ with respect to $v_{i+1}$ have been effectively compensated.

The above predictive compensation method shows that, in general, by use of the m-order synthesis function (200), the predictive compensation of the front $(m-1)$ extracted signals can be obtained at most and the last extracted signal $x_m$ cannot be compensated.

## 7. Numerical simulation and analysis for m-order synthesis function

**Example 7.1.** The given signal is

$$v = v_1 + v_m g_{sm} wgn\left(1, length, -20\right) \tag{209}$$

$$v_1 = v_m \sin\left(\omega t\right) \tag{210}$$

Where $v_m$ and $\omega$ are the amplitude and angular frequency, respectively, $g_{sm}$ is the relative amplitude of Gaussian white noise, $wgn(1, length, -20)$ is a Gaussian white noise function that "1" represents the number of rows, "length" represents the number of columns, and then "$-20$" represents the noise strength of $-20\text{dBW}$.

The other parameters are listed below.

$\omega = 6.28 rad/s$, $v_m = 2$, $h = 0.0005s$

The 3rd-order synthesis function ($m = 3$) is first used to extract the signal and its differential signals.

$$u = fxiao(x_1 - v, x_2, x_3, r, n_1 h) \tag{211}$$

The following algorithm is used to predict and compensate the phase delay and amplitude attenuation.

$$\begin{cases} xiu_1 = x_1 + mn_1 h x_2 + \dfrac{(mn_1 h)^2}{4} x_3 \\ xiu_2 = x_2 + mn_1 h x_3 \end{cases} \tag{212}$$

1) $n_1 = 10$, $g_{sm} = 0.001$

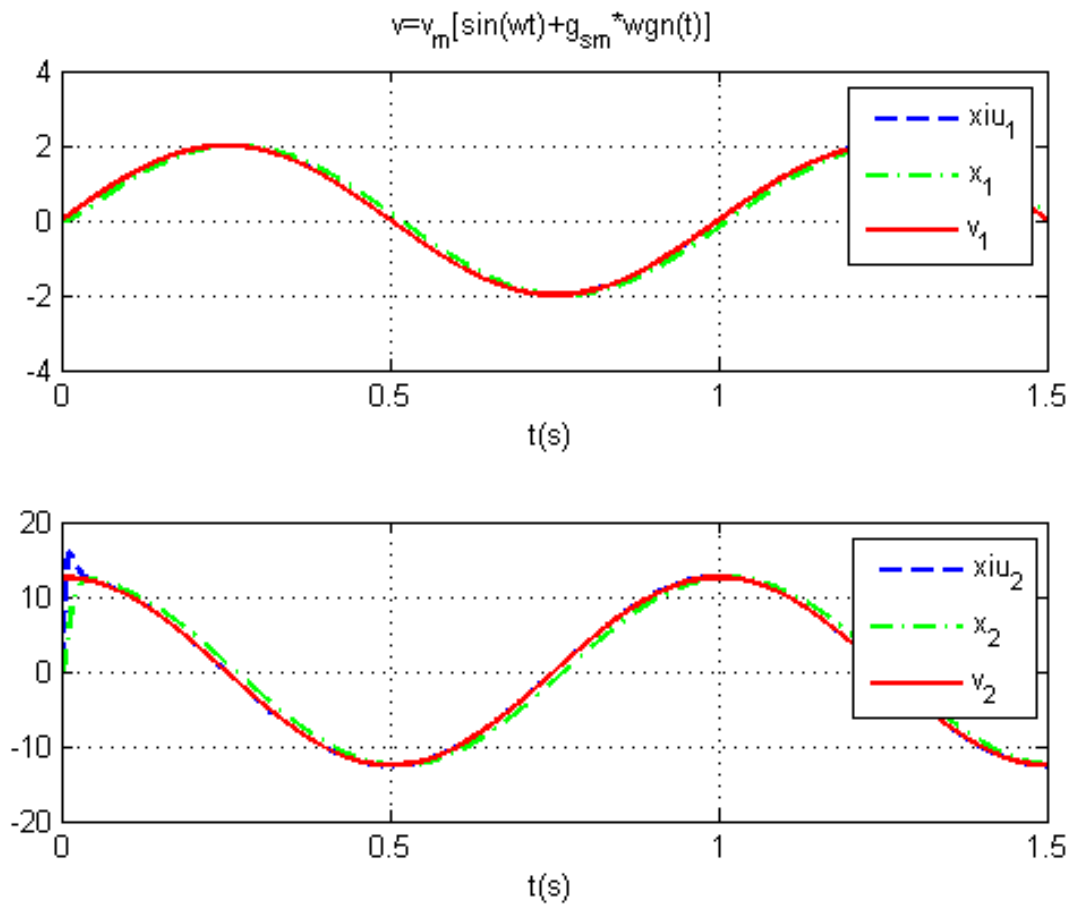

Fig.2. $n_1 = 10$, $g_{sm} = 0.001$ (3rd-order)

In Fig.2, $v_2$ is the differential signal of $v_1$, $x_1$ and $x_2$ are the extracted signal and its differential signal of $v$, and $xiu_1$ and $xiu_2$ are the compensated signals of $x_1$ and $x_2$, respectively. It can be seen that $x_1$ (or $x_2$) is correspondingly behind $v_1$ (or $v_2$) by phase of about $mn_1 h = 30h$, respectively, and after the predictive compensation, the phase and amplitude deviation between $xiu_1$ (or $xiu_2$) and $v_1$ (or $v_2$) are very small in which the phase error is less than one step length $h$ and the relative amplitude error is about $0.1\%$, respectively. This shows that the predictive compensation method (204) is very effective.

2) $n_1 = 10$, $g_{sm} = 0.01$

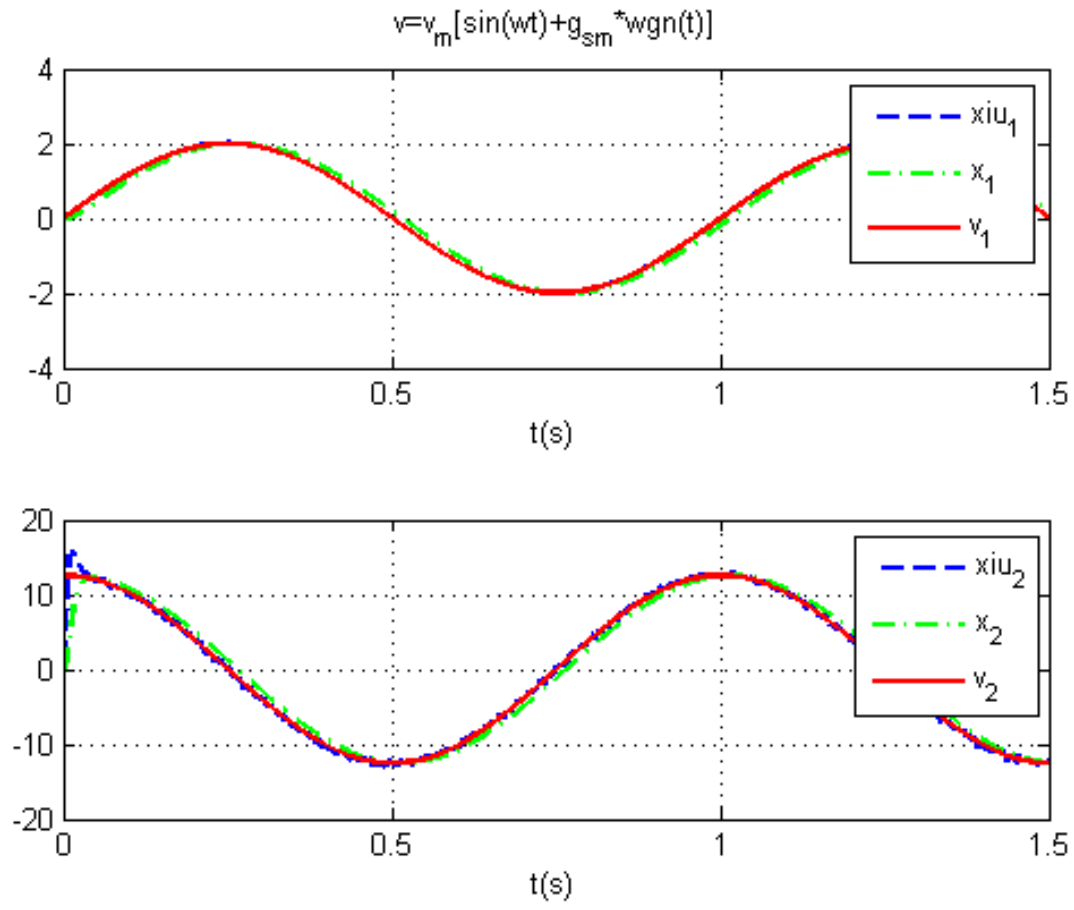

Fig.3. $n_1 = 10$, $g_{sm} = 0.01$ (3rd-order)

3) $n_1 = 10$, $g_{sm} = 0.1$

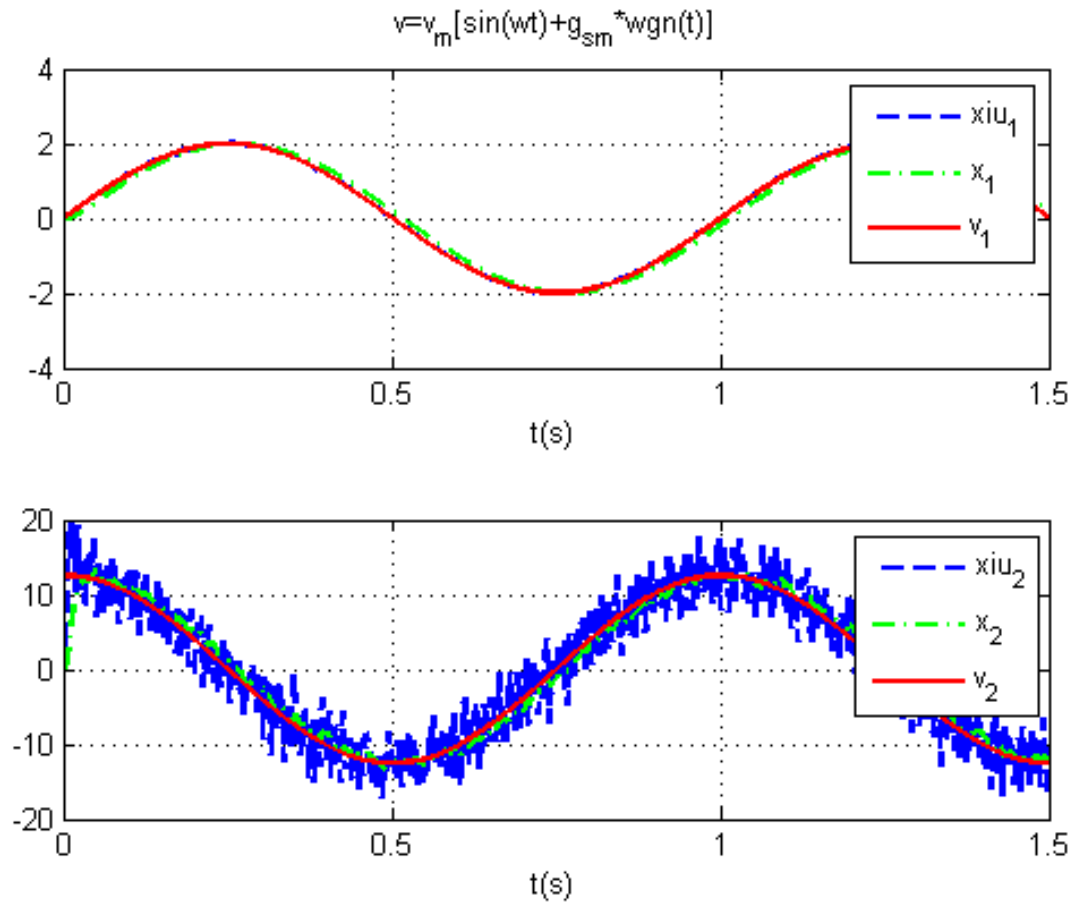

Fig.4. $n_1 = 10$, $g_{sm} = 0.1$ (3rd-order)

When the amplitude of Gaussian white noise is increased, the influence on $x_1$ (or $xiu_1$) is very small; however, the influence on $x_2$ (or $xiu_2$) is greater.

4) $n_1 = 20$, $g_{sm} = 0.01$

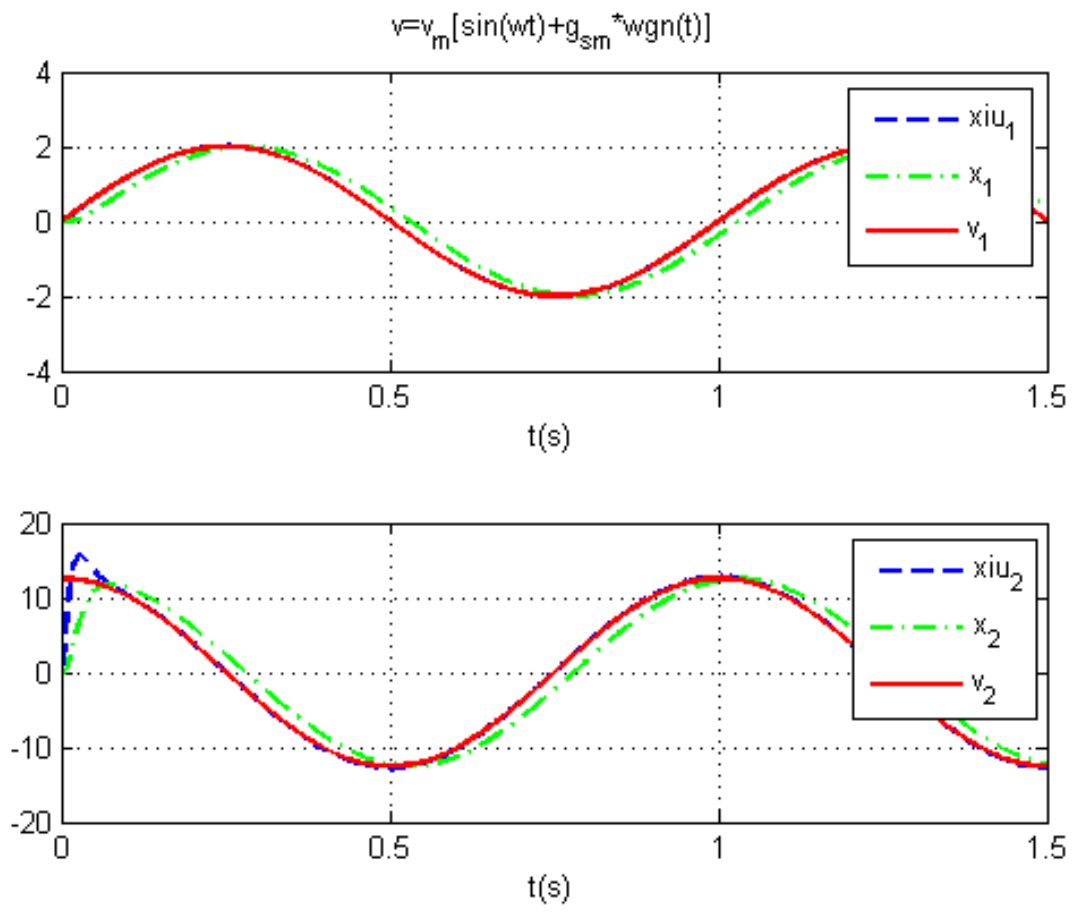

Fig.5. $n_1 = 20$, $g_{sm} = 0.01$ (3rd-order)

5) $n_1 = 30$, $g_{sm} = 0.01$

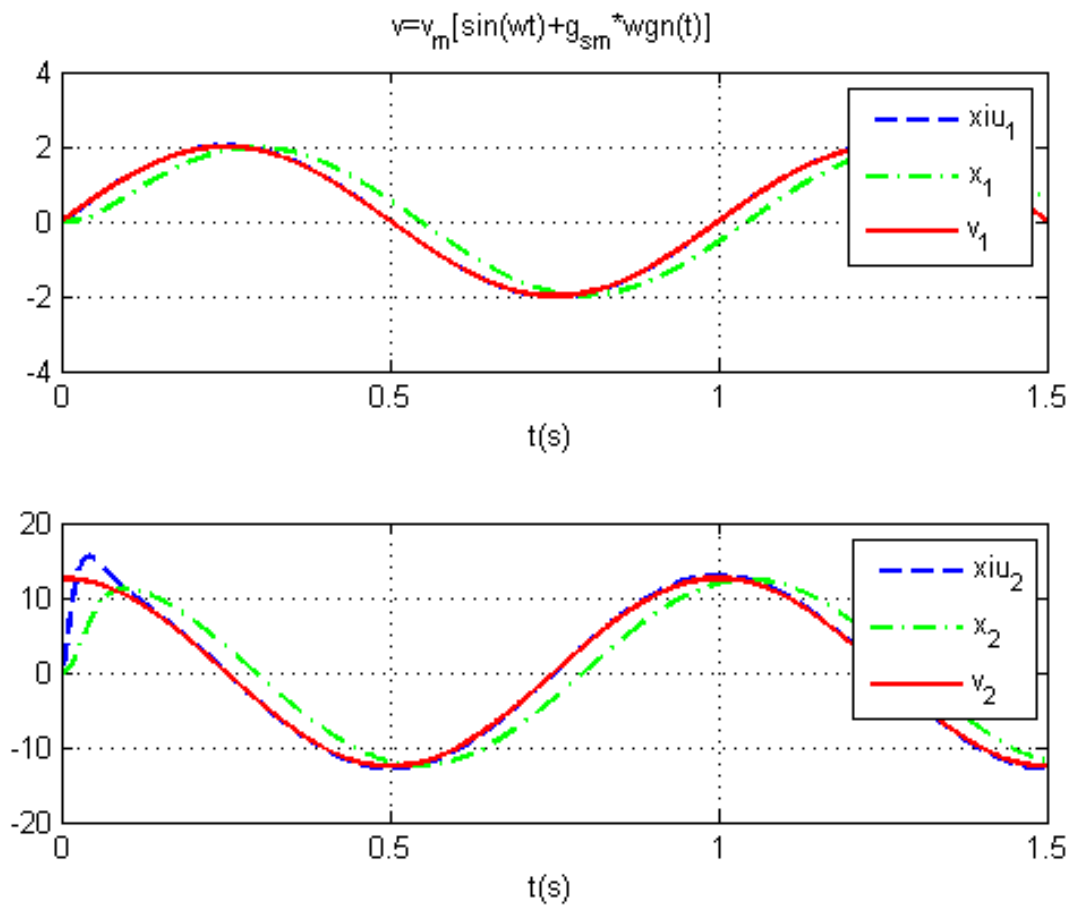

Fig.6. $n_1 = 30$, $g_{sm} = 0.01$ (3rd-order)

6) $n_1 = 40$, $g_{sm} = 0.01$

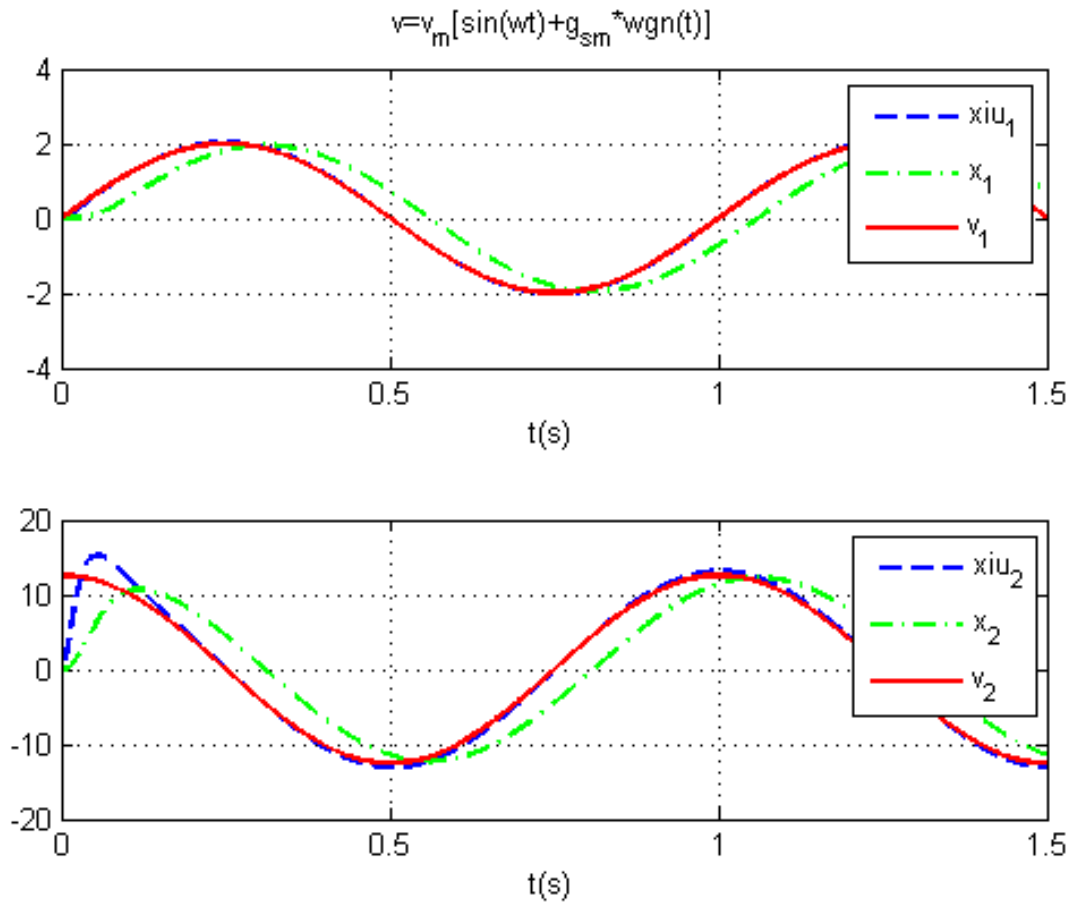

Fig.7. $n_1 = 40$, $g_{sm} = 0.01$ (3rd-order)

It can be seen from Fig.3 and Fig.5-Fig.7 that the phase delay of $x_1$ (or $x_2$) with respect to $v_1$ (or $v_2$) is directly proportional to the filter factor of $n_1$, respectively, and however, the compensated signals $xiu_1$ and $xiu_2$ are substantially independent of the filter factor $n_1$.

**Example 7.2.** The given signal is (209) and (210), the other parameters are the same as in Example 7.1.

The 4th-order synthesis function ($m = 4$) is then used to extract the signal and its differential signals.

$$u = fxiao\left(x_1 - v, x_2, x_3, x_4, r, n_1 h\right) \tag{213}$$

The following algorithm is used to predict and compensate the phase delay and amplitude attenuation.

$$\begin{cases} xiu_1 = x_1 + mn_1 hx_2 + \dfrac{(mn_1 h)^2}{3} x_3 + \left(\dfrac{mn_1 h}{3}\right)^3 x_4 \\ xiu_2 = x_2 + mn_1 hx_3 + \dfrac{(mn_1 h)^2}{4} x_4 \\ xiu_3 = x_3 + mn_1 hx_4 \end{cases} \tag{214}$$

1) $n_1 = 10$, $g_{sm} = 0.001$

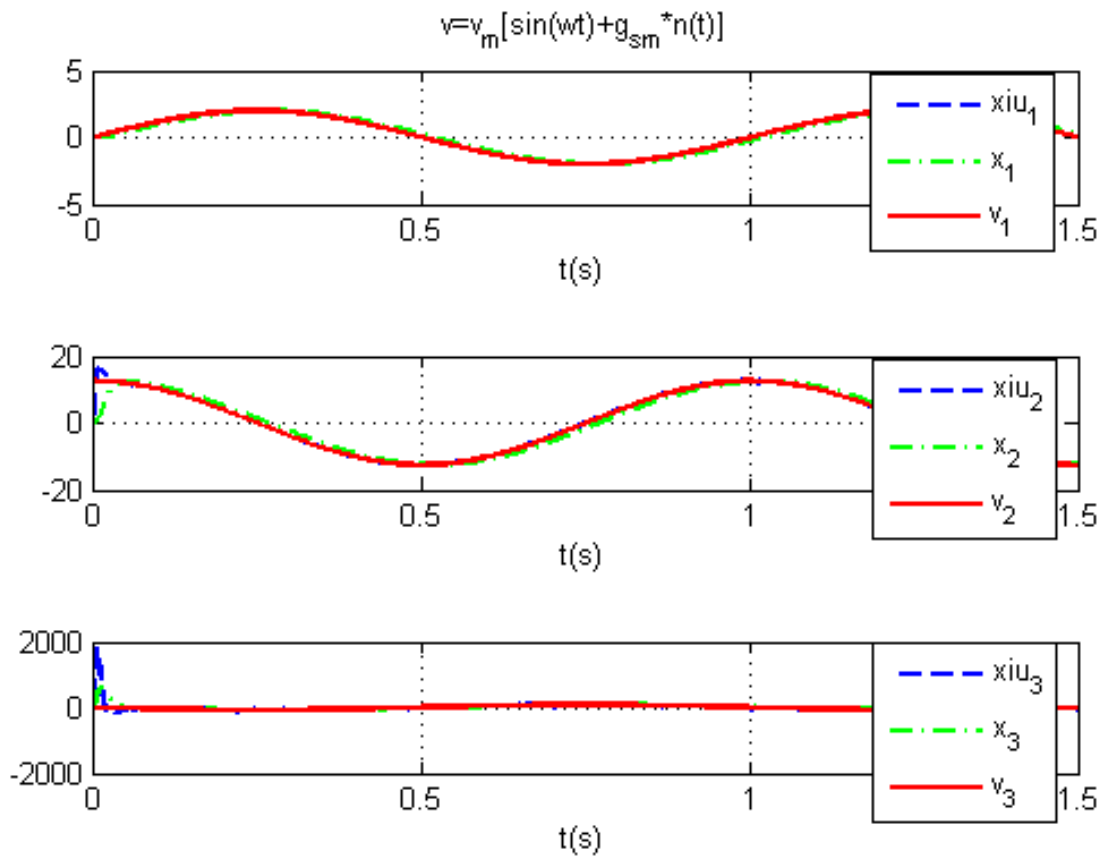

Fig.8. $n_1 = 10$, $g_{sm} = 0.001$ (4th-order)

In Fig.8, $v_2$ and $v_3$ is the 1st-order and 2nd-order differential signals of $v_1$, $x_1$, $x_2$ and $x_3$ are the extracted signal and its differential signals of $v$, and $xiu_1$, $xiu_2$ and $xiu_3$ are the compensated signals of $x_1$, $x_2$ and $x_3$, respectively. It can be seen that $x_1$ (or $x_2, x_3$) is correspondingly behind $v_1$ (or $v_2, v_3$) by phase of about $mn_1 h = 40h$, respectively, and after the predictive compensation, the phase and amplitude deviation between $xiu_1$ (or $xiu_2, xiu_3$) and $v_1$ (or $v_2, v_3$) are very small, respectively. This shows that the predictive compensation method (204) is very effective.

2) $n_1 = 10$, $g_{sm} = 0.01$

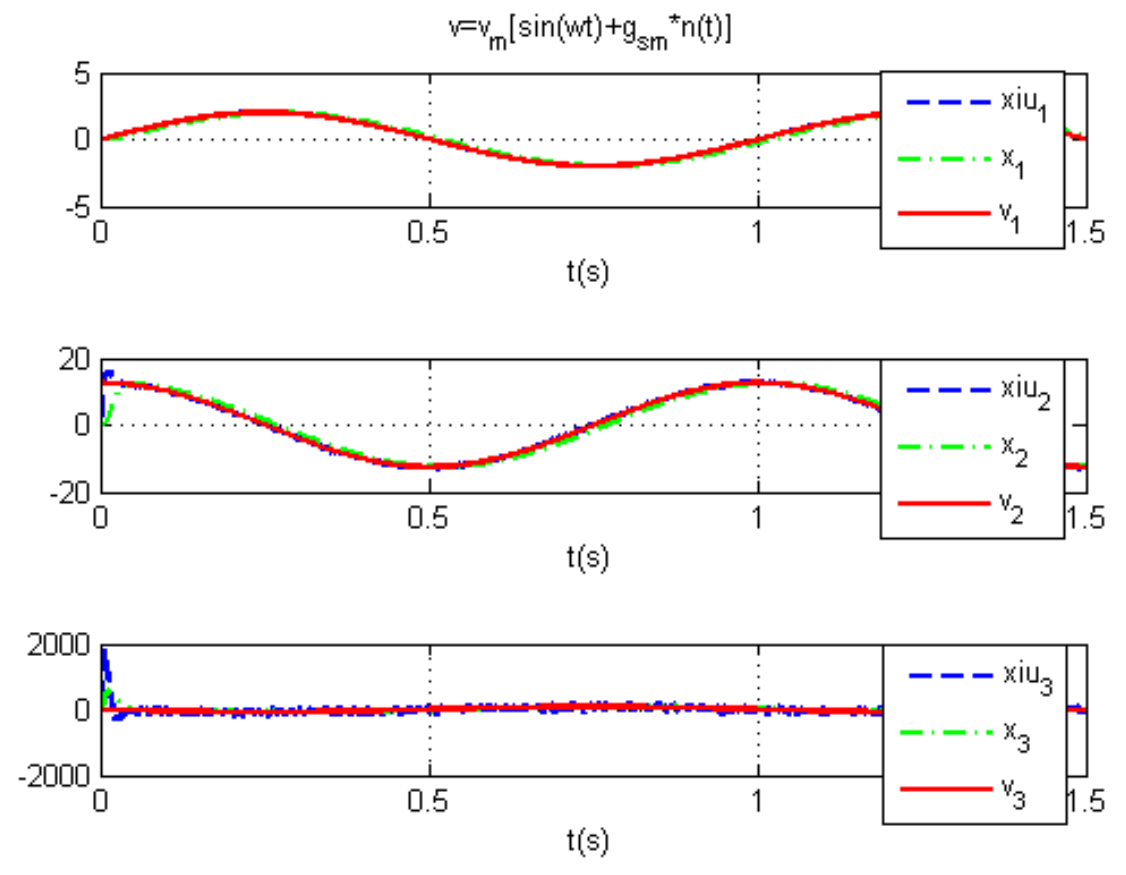

Fig.9. $n_1 = 10$ , $g_{sm} = 0.01$ (4th-order)

3) $n_1 = 10$ , $g_{sm} = 0.1$

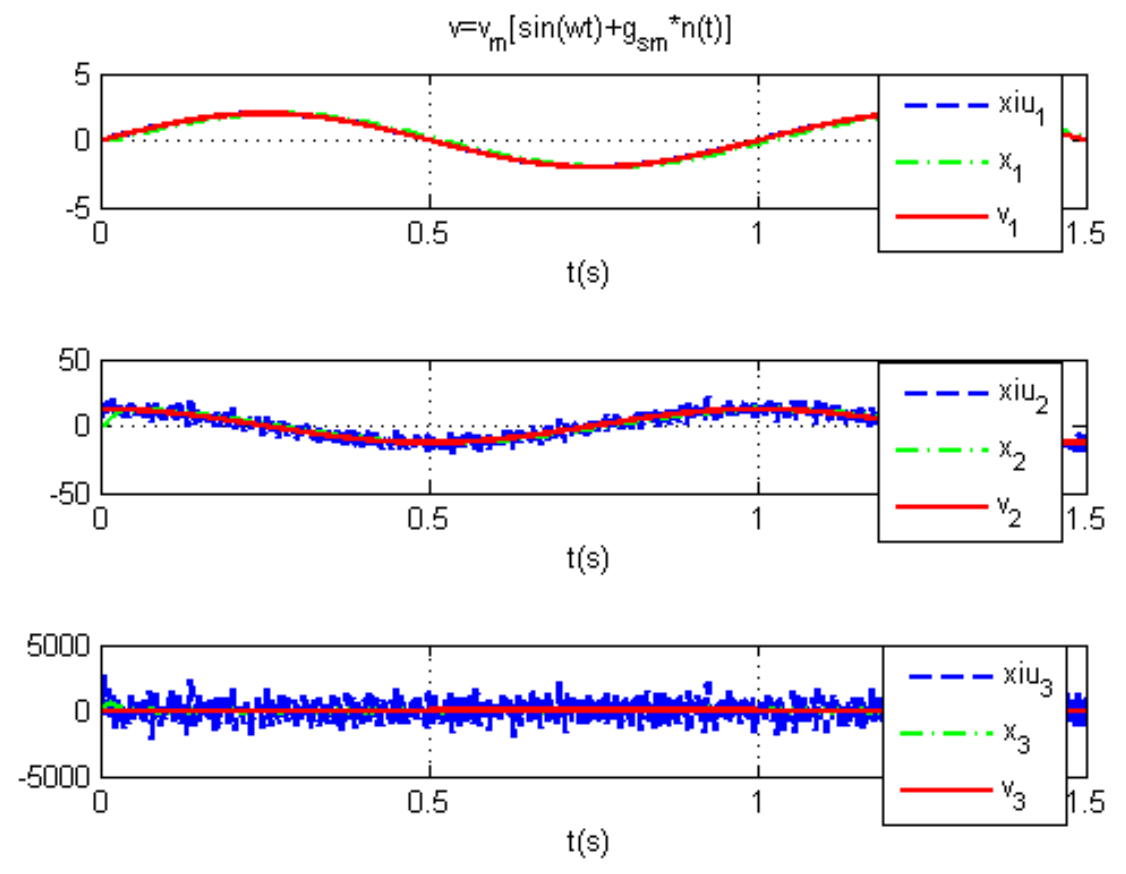

Fig.10. $n_1 = 10$ , $g_{sm} = 0.1$ (4th-order)

When the amplitude of Gaussian white noise is increased, the influence on $x_1$ (or $xiu_1$ ) is very small; however, the influence on $x_2$ (or $xiu_2, x_3, xiu_3$ ) is greater.

4) $n_1 = 20$ , $g_{sm} = 0.01$

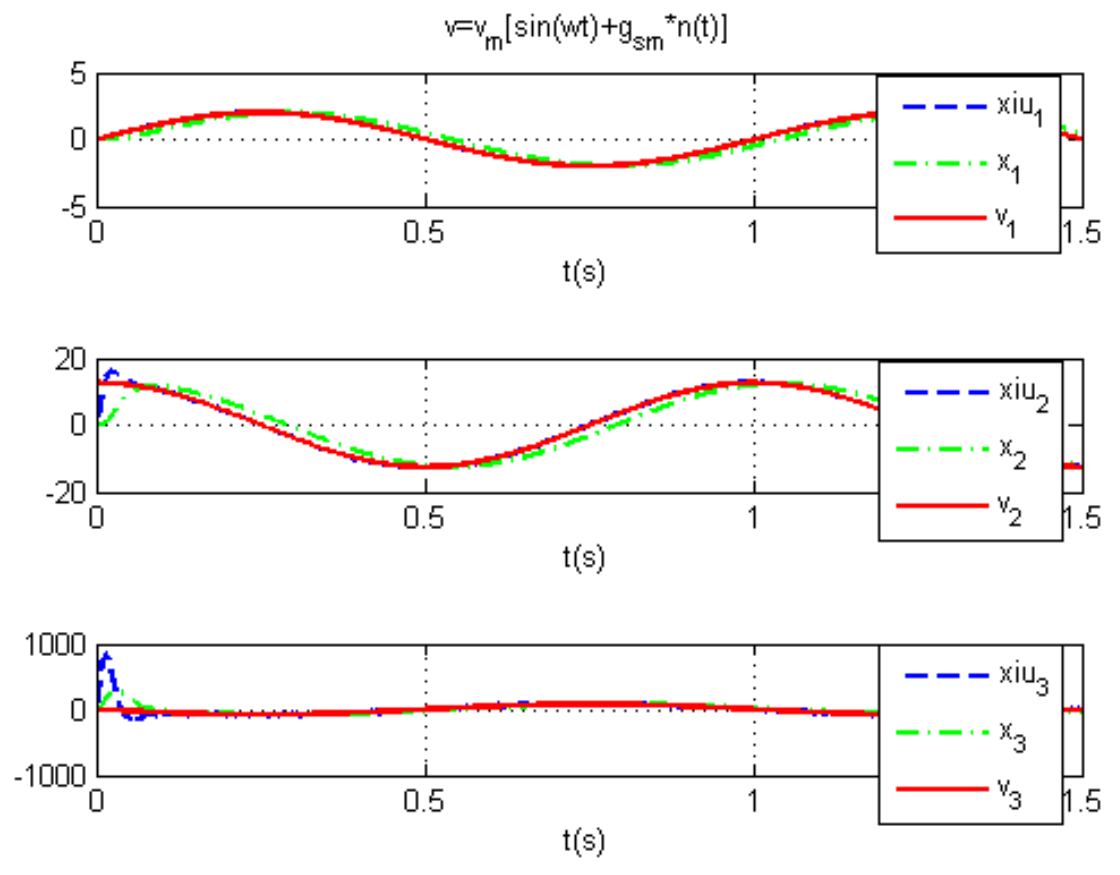

Fig.11. $n_1 = 20$ , $g_{sm} = 0.01$ (4th-order)

5) $n_1 = 30$ , $g_{sm} = 0.01$

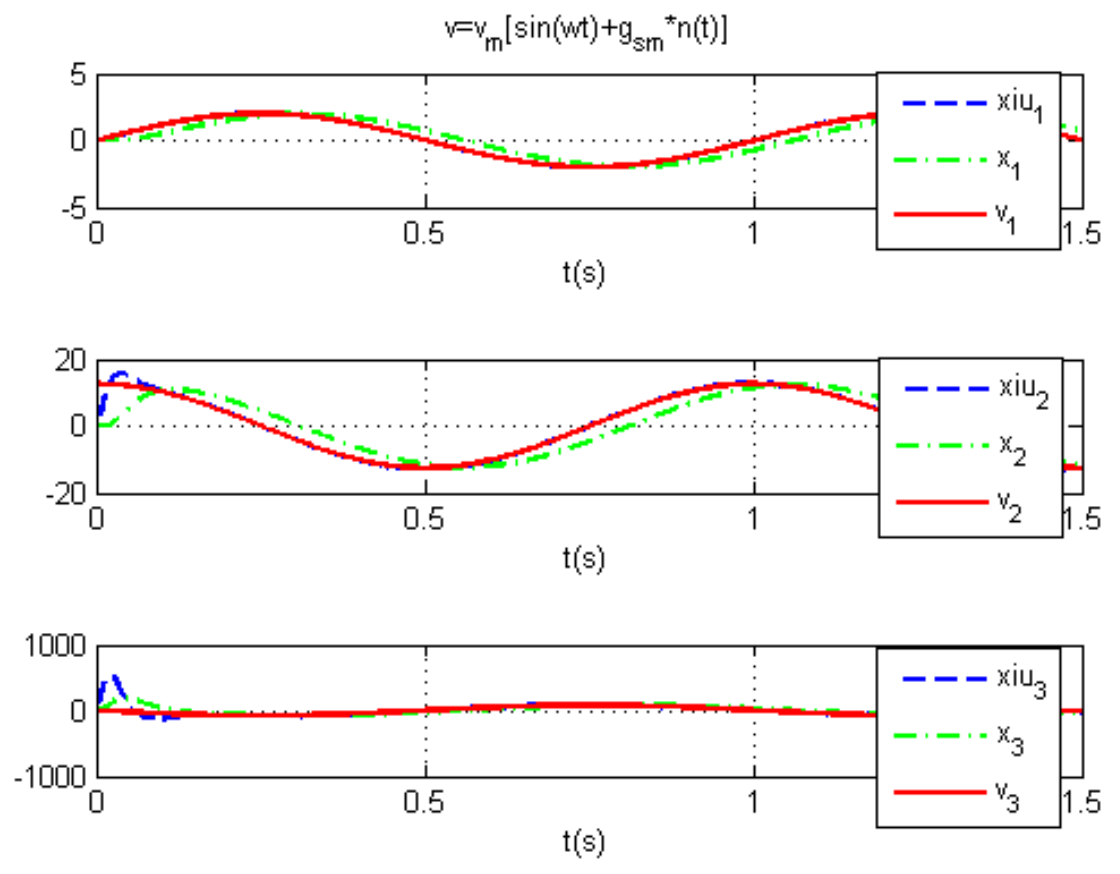

Fig.12. $n_1 = 30$ , $g_{sm} = 0.01$ (4th-order)

6) $n_1 = 40$ , $g_{sm} = 0.01$

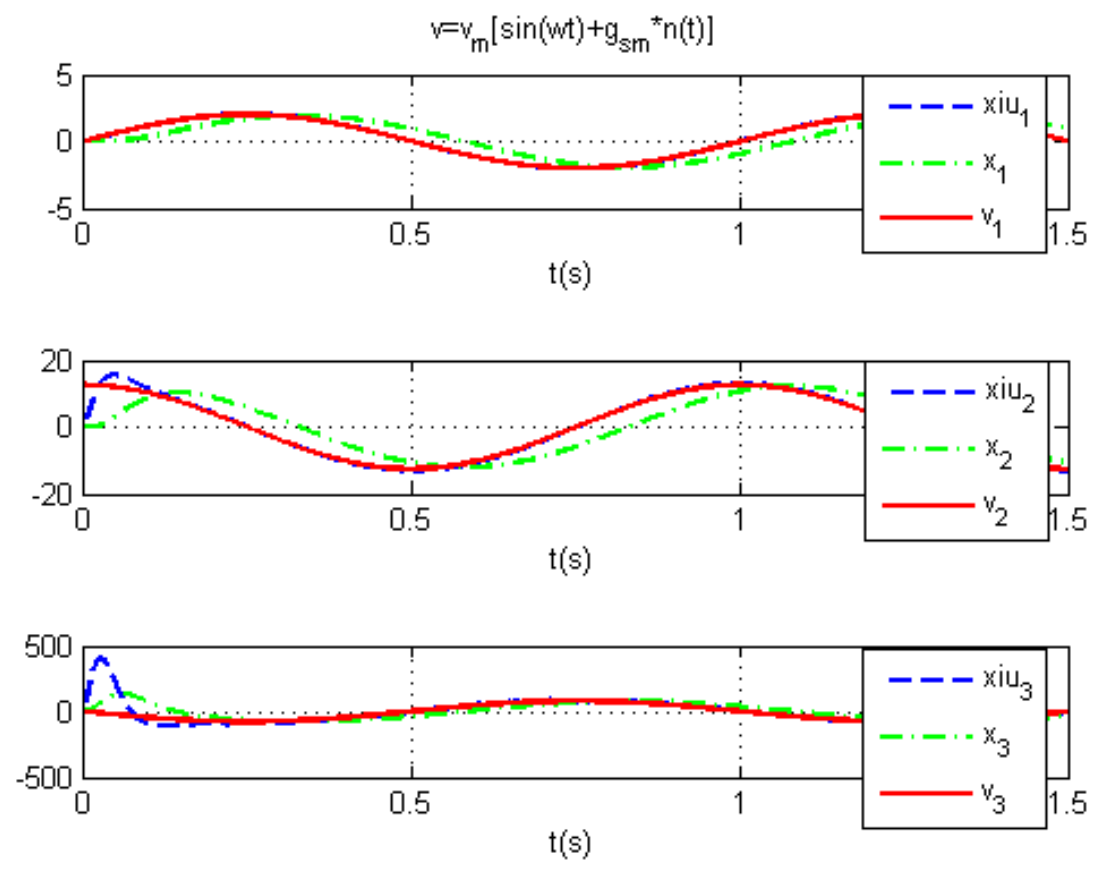

Fig.13. $n_1 = 40$ , $g_{sm} = 0.01$ (4th-order)

It can be seen from Fig.9 and Fig.11-Fig.13 that the phase delay of $x_1$ (or $x_2, x_3$) with respect to $v_1$ (or $v_2, v_3$) is directly proportional to the filter factor of $n_1$, respectively, and however, the compensated signal $xiu_1$ (or $xiu_2, xiu_3$) is substantially independent of the filter factor $n_1$.

## 8. Second-order ADRC technology and analysis

**Definition 8.1.** $m \geq 2$, the m-order integral series system satisfies the following equations.

$$\begin{cases} \dot{x}_1 = x_2 \\ \dot{x}_2 = x_3 \\ \vdots \\ \dot{x}_m = u_0 \\ \chi = x_1 \end{cases} \tag{215}$$

Where, $\chi$ is the output (or feedback) signal, and $u_0$ is the normalized control input.

And the 2nd-order object is

$$\begin{cases} \dot{x}_1 = x_2 \\ \dot{x}_2 = f\left(x_1, x_2, t, w(t)\right) + bu \\ \chi = x_1 \end{cases} \tag{216}$$

Where, $\chi$ is the output signal (or feedback signal), $u$ is the control input, $w(t)$ is the disturbance term of system, $f\left(x_1, x_2, t, w(t)\right)$ is the known (or partially unknown, unknown) total disturbance term (TDT) (no distinction is made between the internal and external disturbances of system (Han, 1998, 2008)), and $b \neq 0$. Set $\xi \triangleq bu$, $u_0 \triangleq f + bu$, $u_0$ is the normalized control input.

The ADRC controller with the estimation and compensation function of TDT is composed of the following parts (Han, 1998, 2008):

1) Transition process arrangement (Han, 1998, 2008)

To arrange the transition process in advance is an effective way to solve the conflict between the overshoot and rapidity of control system, and it is also an effective way to improve the robustness, adaptability and stability of closed-loop system. The specific treatment is that according to the given signal $v$, the transition process $v_1$ and its differential signal $v_2$ are arranged as follows.

$$\begin{cases} u_v = fsun\left(v_1 - v, v_2, r_1, n_1 h\right) \\ \dot{v}_1 = v_2 \\ \dot{v}_2 = u_v \end{cases} \tag{217}$$

Where the parameters of $n_1$ and $r_1$ only depend on the requirements of transition process.

As shown in the previous analysis, the extracted signal $v_1$ and its differential signal $v_2$ also have phase delay and amplitude attenuation.

2) Extended state observer (ESO) (Gao, 2003; Gao, Huang, & Han, 2001; Guo & Zhao, 2016; Han, 1995a, 1998, 2008)

The observed values $y_1$ and $y_2$ for the object states $x_1$ and $x_2$, and the total disturbance term $y_3$ (i.e., extended state) are estimated from the output signal $\chi$ and input signal $u$ of system. The linear ESO is used here to analyze the problem better.

$$\begin{cases} e = y_1 - \chi \\ \dot{y}_1 = y_2 - \beta_{01} e \\ \dot{y}_2 = y_3 - \beta_{02} e + bu \\ \dot{y}_3 = -\beta_{03} e \end{cases} \tag{218}$$

Considering the zero-state response of (218), using the Laplace operator $s$ instead of the differential operator, the transfer relations are obtained as follows.

$$\begin{cases} y_1(s) = \dfrac{\beta_{01}s^2 + \beta_{02}s + \beta_{03}}{\Delta}\chi(s) + \dfrac{s}{\Delta}\xi(s) \\ y_2(s) = \dfrac{\beta_{02}s^2 + \beta_{03}s}{\Delta}\chi(s) + \dfrac{s^2 + \beta_{01}s}{\Delta}\xi(s) \\ y_3(s) = \dfrac{\beta_{03}s^2}{\Delta}\chi(s) - \dfrac{\beta_{03}}{\Delta}\xi(s) \end{cases} \tag{219}$$

Where $\Delta \triangleq s^3 + \beta_{01}s^2 + \beta_{02}s + \beta_{03}$ is the characteristic polynomial of (218), and $\beta_{01}$, $\beta_{02}$, $\beta_{03}$ are a set of parameters.

The characteristic equation of (218) is

$$\Delta = s^3 + \beta_{01}s^2 + \beta_{02}s + \beta_{03} = 0 \tag{220}$$

In order to estimate the object states and TDT, the appropriate parameters of $\beta_{01}$, $\beta_{02}$ and $\beta_{03}$ are need to select, so that (220) is more stable. They can be selected as follows (Gao, Huang, & Han, 2001).

$$\Delta = \left(s + \omega_c\right)^3 = \left(s + \frac{1}{n_2 h}\right)^3 \tag{221}$$

$$\begin{cases} \beta_{01} = 3\omega_c = \dfrac{3}{n_2 h} \\ \beta_{02} = 3\omega_c^2 = \dfrac{3}{\left(n_2 h\right)^2} \\ \beta_{03} = \omega_c^3 = \dfrac{1}{\left(n_2 h\right)^3} \end{cases} \tag{222}$$

Substitute (222) into (219) and obtain

$$\begin{cases} y_1(s) = \dfrac{1+3n_2hs+3(n_2h)^2 s^2}{(1+n_2hs)^3}\chi(s) + \dfrac{(n_2h)^3 s}{(1+n_2hs)^3}\xi(s) \\ y_2(s) = \dfrac{s+3n_2hs^2}{(1+n_2hs)^3}\chi(s) + \dfrac{(n_2h)^3 s^2 + 3(n_2h)^2 s}{(1+n_2hs)^3}\xi(s) \\ y_3(s) = \dfrac{s^2}{(1+n_2hs)^3}\chi(s) - \dfrac{1}{(1+n_2hs)^3}\xi(s) \end{cases} \quad (223)$$

When $h$ is small enough, the 1st-order approximation with respect to $h$ can be done on the right of (223) as follows.

$$\begin{cases} y_1(s) \approx \dfrac{1+3n_2hs}{1+3n_2hs}x_1(s) = x_1(s) \\ y_2(s) \approx \dfrac{1+3n_2hs}{1+3n_2hs}x_2(s) = x_2(s) \\ y_3(s) \approx \dfrac{s^2}{1+3n_2hs}x_1(s) - \dfrac{1}{1+3n_2hs}\xi(s) \end{cases} \quad (224)$$

It can be seen from (224) that $y_1$ and $y_2$ are the approximate observations of object states $x_1$ and $x_2$, with little phase delay and amplitude attenuation; and $y_3$ is the filtered estimate of TDT $f$, with the phase delay and amplitude attenuation. Moreover, the extracted signals obtained in (218) cannot be used to predict and compensate the signal $y_3$.

3) Nonlinear error feedback control law (Han, 1995b, 1998, 2008)

The system's state errors are

$$\begin{cases} e_1 = v_1 - y_1 \\ e_2 = v_2 - y_2 \end{cases} \quad (225)$$

The nonlinear error feedback control law is based on the errors $e_1$ and $e_2$ to determine the normalized control input $u_0$ of 2nd-order integral series system, which can adopt nonlinear PID, or $fhan(e_1, e_2, r_3, n_3h)$ (Han, 2008; Han & Yuan, 1999), etc. In order to facilitate the analysis, the control input $u_0$ is obtained by use of the 2nd-order error-based linear synthesis function.

$$u_0 = \frac{e_1 + 2n_3he_2}{(n_3h)^2} \quad (226)$$

As a result, the 2nd-order linear time optimal control system is obtained as follows.

$$\begin{cases} u_0 = \dfrac{(v_1 - x_1) + 2n_3h(\dot{v}_1 - \dot{x}_1)}{(n_3h)^2} \\ \dot{x}_1 = x_2 \\ \dot{x}_2 = u_0 \end{cases} \quad (227)$$

Considering the zero-state response of (227), using the Laplace operator $s$ instead of the differential operator, the transfer function is obtained as follows.

$$\frac{x_1(s)}{v_1(s)} = \frac{1+2n_3hs}{(1+n_3hs)^2} \quad (228)$$

(228) shows that using the linear error feedback law will make the relative order of closed-loop system become one by introducing an additional zero point, which the open-loop system has the relative order of two. The additional closed-loop zero point will affect the system's dynamic response process, which is why the approximate 2nd-order synthesis function *fhan*() requires an additional parameter $c$ (control input is $fhan(e_1, ce_2, r_3, n_3h)$) to adjust the dynamic response process in practical application (Han, 2008). Similarly, using the m-order error-based linear synthesis function in m-order linear system (215) will introduce $(m-1)$ additional zero points in closed-loop system (whose relative order is also one), which have an impact on the system's response process. This should attract enough attention.

When $h$ is small enough, the 1st-order approximation with respect to $h$ can be done on the right of (228) as follows.

$$\frac{x_1(s)}{v_1(s)} \approx \frac{1+2n_3hs}{1+2n_3hs} = 1 \quad (229)$$

(229) shows that there is a feature which is not easy to notice after using the linear error feedback law, i.e., the object state $x_1$ has almost no phase delay and amplitude attenuation for the arranged transition process $v_1$; however, the relative order of state $x_2$ relative to $v_1$ is zero, which shows that $x_2$ does not have the ability to suppress noise in $v_1$. Similarly, the state $x_1$ of closed-loop system has almost no phase delay and amplitude attenuation for the arranged transition process $v_1$ after using the m-order error-based linear synthesis function (198) in m-order linear system (215); however, the relative order of state $x_{i+1}(1 \le i \le m-1)$ relative to $v_1$ is not a positive integer, which also shows that $x_{i+1}$ does not have the ability to suppress noise in $v_1$.

4) Compensation of TDT (Han, 1995b, 1998, 2008)

For the normalized control input $u_0$, the actual control input $u$ of system is determined by the compensation of TDT's estimate $y_3$.

$$u = \frac{u_0 - y_3}{b} \tag{230}$$

Using the TDT's estimate to compensate the control input is one of most groundbreaking ideas in active disturbance rejection controller. Active disturbance rejection characteristics are referred to the real-time estimation and compensation functions of total disturbances, that is, anti-interference functions. This idea deserves to be inherited and carried forward.

## 9. New m-order ADRC theory

The m-order object is

$$\begin{cases} \dot{x}_1 = x_2 \\ \dot{x}_2 = x_3 \\ \vdots \\ \dot{x}_m = f\left(x_1, x_2, \cdots, x_m, t, w(t)\right) + bu \\ \chi = x_1 \end{cases} \tag{231}$$

Where, $\chi$ is the output signal (or feedback signal), $u$ is the control input, $w(t)$ is the disturbance term of system, $f\left(x_1, x_2, \cdots, x_m, t, w(t)\right)$ is the known (or partially unknown, unknown) TDT (no distinction is made between the internal and external disturbances of system (Han, 1998, 2008)), and $b \neq 0$ . Set $\xi \triangleq bu$ , $u_0 \triangleq f + bu$ , $u_0$ is the normalized control input.

For system (231), a new m-order discrete-time ADRCT is proposed in the following.

### 9.1. *Transition process arrangement*

According to the given signal $v$, the transition process $v_1(k)$ can be arranged by use of the 2nd-order synthesis function (for the specific formula, see (217)) or the $m_1(\geq 2)$-order synthesis function. Mentioned earlier, the parameters of $m_1$, $n_1$, and $r_1$ only depend on the requirements of transition process, which shows that the design of transition process has nothing to do with the ESO and controller.

In particular, the arranged transition process can also be any given curve based on the characteristics of object, not just the methods in this section.

### 9.2. *Improved extended state observer*

In order to predict and compensate the TDT's estimate when needed, the $(m+2)$-order synthesis function is used to extract the signal $y_1(k)$ and its differential signals $y_{i+1}(k)(1 \leq i \leq m+1)$ for output $\chi$ at $k$ moment as follows.

$$\begin{cases} u_y(k-1) = \\ fxiao\left(y_1(k-1) - \chi, y_2(k-1), \cdots y_{m+2}(k-1), r_2, n_2 h\right) \\ y_1(k) = y_1(k-1) + hy_2(k-1) \\ y_2(k) = y_2(k-1) + hy_3(k-1) \\ \vdots \\ y_{m+2}(k) = y_{m+2}(k-1) + hu_y(k-1) \end{cases} \tag{232}$$

The extracted signal $y_{i+1}(k)(0 \leq i \leq m-1)$ used to calculate the normalized control input can be compensated for its phase delay when needed (see (204)), as follows.

$$\hat{y}_{i+1}(k) = \sum_{\mu=0}^{m+1-i} C_{m+1-i}^{\mu} \left[\frac{(m+2)n_2 h}{m+1-i}\right]^{\mu} y_{i+1+\mu}(k) \tag{233}$$

In addition, in order to obtain the estimate of TDT, the $(m+2)$-order synthesis function with the same parameters as (233) must also be used to extract the signal $\xi_1(k-1)$ for the signal $\xi$ at $(k-1)$ moment (at this time, there is no signal $\xi$ at $k$ moment).

$$\begin{cases} u_\xi(k-2) = \\ fxiao\left(\xi_1(k-2) - \xi, \xi_2(k-2), \cdots \xi_{m+2}(k-2), r_2, n_2 h\right) \\ \xi_1(k-1) = \xi_1(k-2) + h\xi_2(k-2) \\ \xi_2(k-1) = \xi_2(k-2) + h\xi_3(k-2) \\ \vdots \\ \xi_{m+2}(k-1) = \xi_{m+2}(k-2) + hu_\xi(k-2) \end{cases} \tag{234}$$

If necessary, the signal $\xi_1(k-1)$ can be compensated for its phase delay as follows.

$$\hat{\xi}_1(k) = \sum_{\mu=0}^{m+1} C_{m+1}^{\mu} \left[\frac{h + (m+2)n_2 h}{m+1}\right]^{\mu} \xi_{\mu+1}(k-1) \tag{235}$$

The signal $\xi_1(k-1)$ is compensated for its one-step phase delay of $h$ as follows.

$$\bar{\xi}_1(k) = \sum_{\mu=0}^{m+1} C_{m+1}^{\mu} \left[\frac{h}{m+1}\right]^{\mu} \xi_{\mu+1}(k-1) \tag{236}$$

Correspondingly, $y_{m+1}(k-1)$ or $y_{m+1}(k)$ can be compensated for their phase delays if necessary.

$$\hat{y}_{m+1}(k) = y_{m+1}(k-1) + \left[(m+2)n_2 + 1\right] hy_{m+2}(k-1) \tag{237}$$

$$\hat{y}_{m+1}(k) = y_{m+1}(k) + (m+2)n_2 hy_{m+2}(k) \tag{238}$$

The total disturbance term $\bar{f}(k)$ can be indirectly observed as follows.

$$\bar{f}(k)=y_{m+1}(k)-\bar{\xi}_1(k) \tag{239}$$

Sometimes the approximate formula is used in the following.

$$\bar{f}(k)\approx y_{m+1}(k-1)-\xi_1(k-1) \tag{240}$$

$$\bar{f}(k)\approx y_{m+1}(k)-\xi_1(k-1) \tag{241}$$

$\bar{f}(k)$ in (240) is behind of $f(k)$ by phase of about $\left[(m+2)n_2+1\right]h$, and if necessary, can be compensated for its phase delay as follows.

$$\hat{f}(k)=\hat{y}_{m+1}(k)-\hat{\xi}_1(k) \tag{242}$$

Since $\hat{y}_{m+1}(k)$ and $\hat{\xi}_1(k)$ are the approximate values of $x_1^{(m)}(k)$ and $\xi$ at $k$ moment, respectively, the prediction signal $\hat{f}(k)$ has substantially no phase delay and amplitude attenuation for $f(k)$.

$$\hat{f}(k)\approx f(k) \tag{243}$$

In the actual discrete computation and control, the prediction signal $\hat{f}(k)$ with fast response in the dynamic process will cause the system to oscillate violently, and on the contrary, its observation accuracy will be much lower than $\bar{f}(k)$. Thus, instead of (242), this paper suggests using (240) or (241), of which (241) is the most concise and practical. (242) and (243) are only used to proof the corresponding theorem in section 9.5.

In this section, we can consider $y_{m+1}(k)$ as the corresponding extended state, and $\xi_1(k)$ is the filtered extracted signal that is forced to be filtered for the indirect estimation of TDT. The TDT is always obtained by indirect observation.

In (232), the order of synthesis function is set to $m_2$. The parameters of $m_2$, $n_2$ and $r_2$ only depend on the requirements of suppressing the noise in feedback signal during the process of extracting the signal and its differential signals, and have nothing to do with the design of controller. However, the same parameters selected in (234) as (232) is only used to estimate the TDT. This shows that the design of ESO has nothing to do with the arranged transition process and controller.

### 9.3. *Feedback control law*

The m-order tracking-form synthesis function is used to compute the normalized control input $u_0$ as follows.

$$u_0(k)=fxiao\left(y_1(k)-v_1(k),y_2(k),\cdots,y_m(k),r_3,n_3h\right) \tag{244}$$

The signal $y_1(k)$ is behind of $v_1(k)$ by phase of about $mn_3h$, and if necessary, can be compensated for its phase delay as follows.

$$u_0(k)=fxiao\left(\hat{y}_1(k)-v_1(k),\hat{y}_2(k),\cdots,\hat{y}_m(k),r_3,n_3h\right) \tag{245}$$

The signal $\hat{y}_1(k)$ is also behind of $v_1(k)$ by phase of about $mn_3h$. There are two ways to deal with this, either one of them can effectively compensate for phase delay and amplitude attenuation.

1) The predictive signal $\bar{v}_1(k)$ is obtained by compensating the signal $v_1(k)$.

$$\bar{v}_1(k)=v_1(k)+(mn_3h)v_2(k) \tag{246}$$

$$u_0(k)=fxiao\left(\hat{y}_1(k)-\bar{v}_1(k),\hat{y}_2(k),\cdots,\hat{y}_m(k),r_3,n_3h\right) \tag{247}$$

2) The predictive signal $\bar{y}_{i+1}(k)$ is obtained by recompensating the signal $\hat{y}_{i+1}(k)$.

$$\bar{y}_{i+1}(k)=\sum_{\mu=0}^{m+1-i}C_{m+1-i}^{\mu}\left[\frac{(m+2)n_2h+mn_3h}{m+1-i}\right]^{\mu}y_{i+1+\mu}(k) \tag{248}$$

$$u_0(k)=fxiao\left(\bar{y}_1(k)-v_1(k),\bar{y}_2(k),\cdots,\bar{y}_m(k),r_3,n_3h\right) \tag{249}$$

At this time, $\hat{y}_1(k)$ (which is an approximation of $x_1(k)$) has substantially no phase delay and amplitude attenuation for the signal $v_1(k)$.

Relatively, method 1) (i.e., (246) and (247)) is simple and easy to implement, and it is recommended to select first.

In (244), the order of synthesis function is $m_3=m$. The parameters of $m_3$, $n_3$ and $r_3$ only depend on the dynamic response performance of system, which indicates that the design of feedback controller have nothing to do with the arranged transition process and ESO.

### 9.4. *Compensation of TDT*

For the normalized control input $u_0(k)$, the actual control input $u(k)$ is determined by the compensation of TDT's estimate $\bar{f}(k)$ (Han, 2008).

$$u(k)=\frac{u_0(k)-\bar{f}(k)}{b(k)} \tag{250}$$

The actual control input $u(k)$ is obtained as follows if the predictive signal $\hat{f}(k)$ is used to compensate.

$$u(k)=\frac{u_0(k)-\hat{f}(k)}{b(k)} \tag{251}$$

In actual discrete computation and control, it is recommended to select (250) instead of (251) because the algorithm is concise and practical.

### 9.5. Dynamic compensation linearization method and its proof

When the TDT $f(k)$ is known, the actual control input of system is taken as follows.

$$u(k)=\frac{u_0(k)-f(k)}{b(k)} \tag{252}$$

The system (231) can be compensated into the m-order integral series system (215) by use of (252). This process is called the direct feedback linearization method (Han, 1981, 2008).

In general, the TDT is partially unknown or unknown, and then the process of compensating the system (231) into the integral series system (215) through (251) is called the dynamic compensation linearization method (Han, 2008). The correctness of dynamic compensation linearization method is proved in the following.

For the $k$ moment, using (242) and (243), there is

$$x_1^{(m)}(k)=f(k)+\xi(k)=f(k)+\left[u_0(k)-\hat{f}(k)\right]$$

$$\approx f(k)+\left[u_0(k)-f(k)\right]=u_0(k) \tag{253}$$

As a result, it is an integral series system (215), which proves the dynamic compensation linearization method.

### 9.6. Summary

Based on the m-order tracking-form synthesis function, the new structures of ESO and feedback controller satisfying the separation principle are reconstructed to construct a new ADRCT. The above derivation and analysis show that the new ADRCT has the following characteristics:

1) The method by estimating and compensating the TDT to counteract the influence of disturbance will be expected to have good robustness, immunity and adaptability.

2) The estimate of TDT has phase delay and amplitude attenuation.

3) The structure of controller inherits the distinctive characteristics of m-order tracking-form synthesis function, and in the linear case, the m-order open-loop system (the relative order of m) can be controlled to one m-order closed-loop system without the additional closed-loop zero point introduced.

4) The parameters of controller and observer are completely independent to each other, that is, the separation principle in modern control theory can be completely satisfied, so that the controller and observer can be independently designed.

## 10. Numerical simulation and analysis for ADRCT

**Example 10.1.** The 3rd-order controlled object is

$$\begin{cases} Dx_1=x_2 \\ Dx_2=x_3 \\ Dx_3=f(t)+u \\ \chi=x_1+v_m g_{sm1} wgn(1,length,-20) \end{cases} \tag{254}$$

$$f(t)=f_1(t)+20g_{sm} wgn(1,length,-20) \tag{255}$$

$$f_1(t)=20[\sin(0.2\omega t)+\sin(0.1\omega t)+\sin(0.05\omega t)] \tag{256}$$

The given signal is $v$, and specific parameters are as follows.

$v=v_m=2$ , $\omega=6.28\text{rad}/\text{s}$ , $h=0.0005s$ , $g_{sm}=g_{sm1}$ ; arranged transition process: $m_1=5$ , $n_1=20$ ; ESO : $m_2=5$ , $n_2=30$ ; controller: $m_3=3$ , $n_3=500$

1) $g_{sm}=0.001$

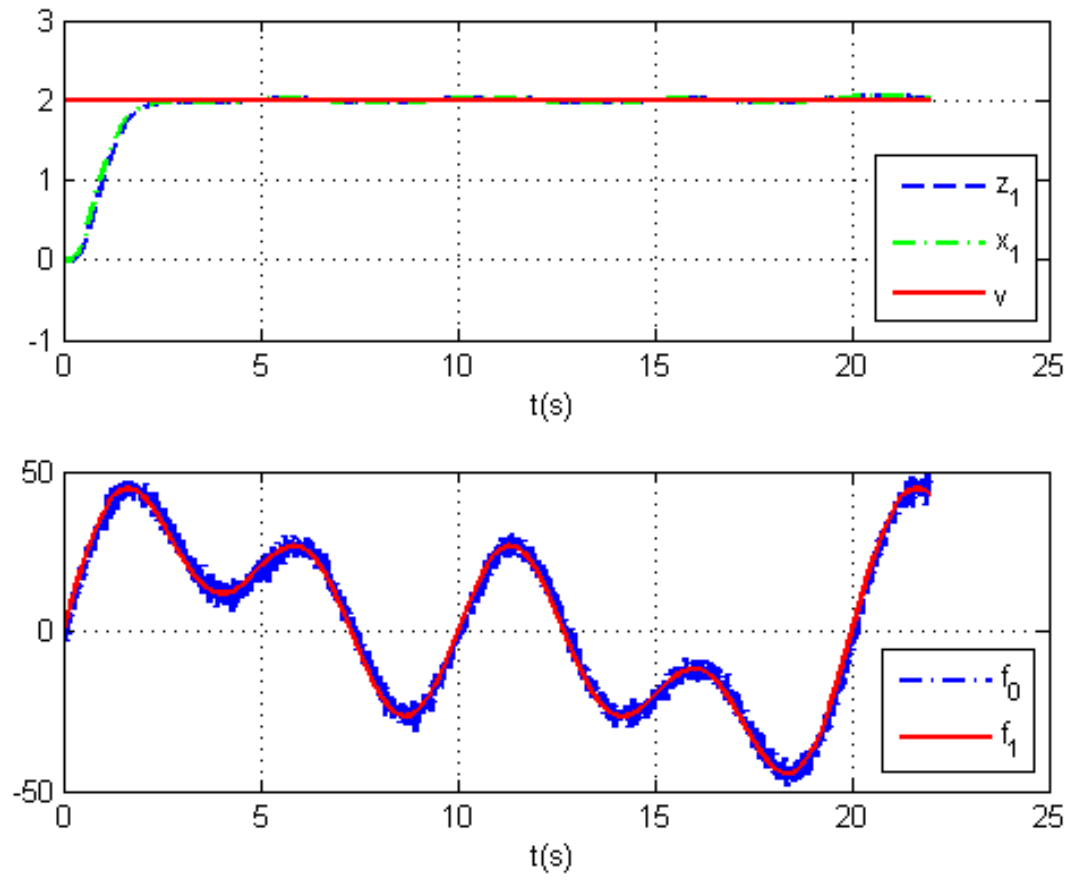

Fig.14. Step response ( $g_{sm}=0.001$ )

2) $g_{sm}=0.01$

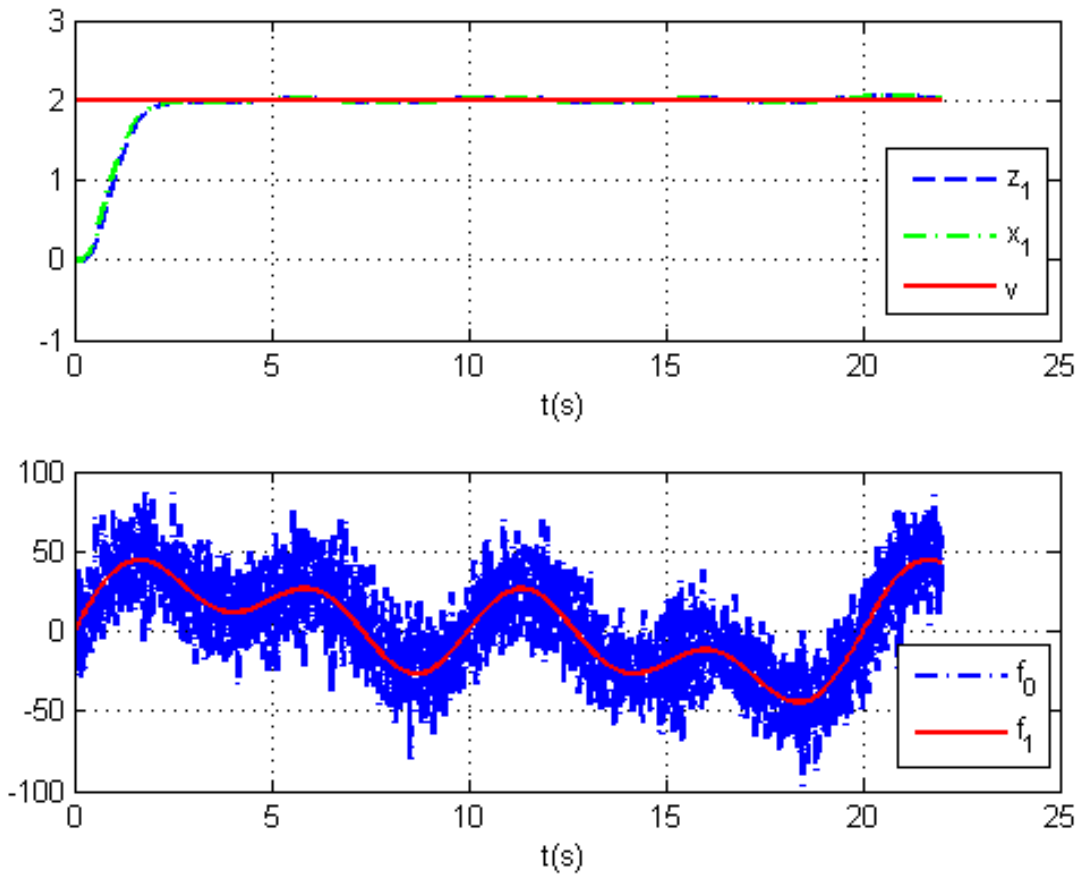

Fig.15. Step response ( $g_{sm} = 0.01$ )

3) $g_{sm} = 0.1$

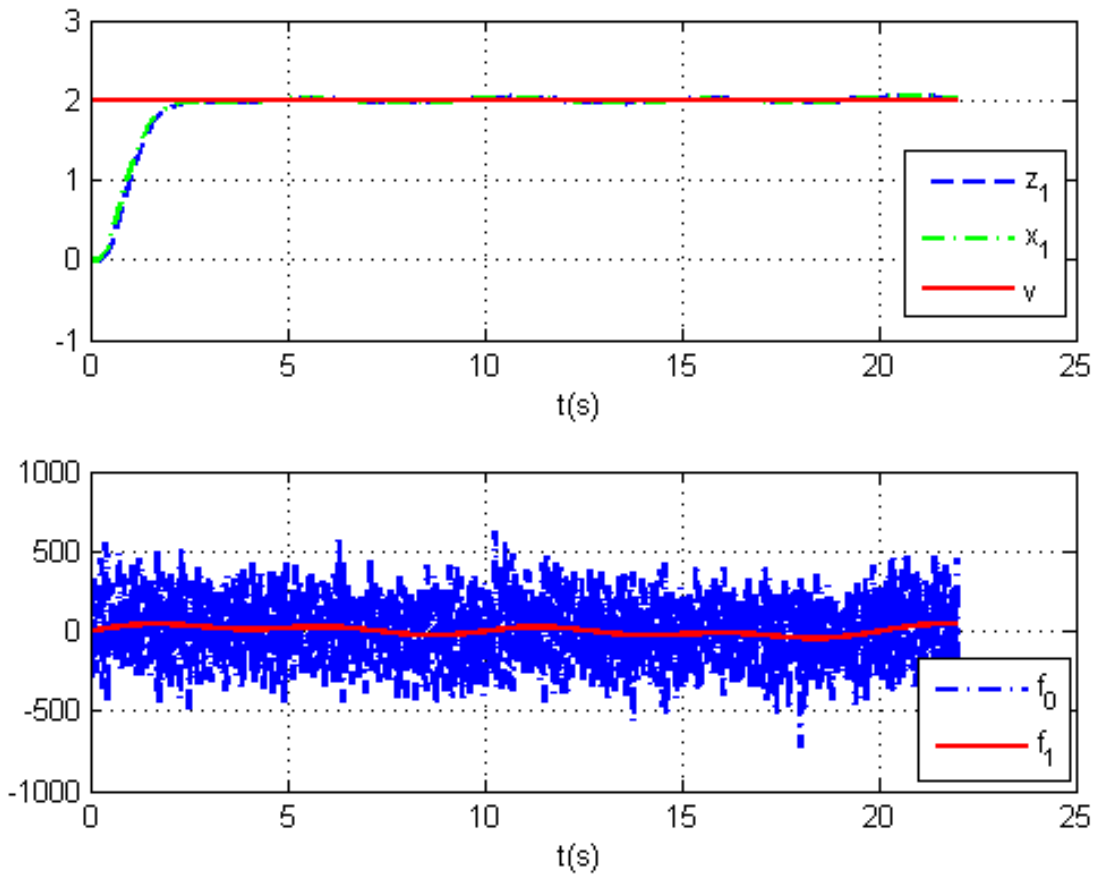

Fig.16. Step response ( $g_{sm} = 0.1$ )

In Fig.14-Fig.16, $z_1$ is the extracted signal for feedback signal $\chi$, $f_1$ represents $f_1(t)$, and $f_0$ is the indirect estimate of TDT $f(t)$. It can be seen that the TDT's estimation algorithm (241) can well counteract the influence of disturbance, and has a strong capability to suppress noise.

**Example 10.2.** The 3rd-order controlled object is

$$\begin{cases} Dx_1 = x_2 \\ Dx_2 = x_3 \\ Dx_3 = f(t) + u \\ \chi = x_1 + v_m g_{sm1} wgn(1, length, -20) \end{cases} \tag{257}$$

$$f(t) = f_1(t) + 20 g_{sm} wgn(1, length, -20) \tag{258}$$

$$f_1(t) = 20 sign(\sin(0.2\omega t)) +$$

$$20 sign(\sin(0.1\omega t)) + 20 sign(\sin(0.05\omega t)) \tag{259}$$

The given signal is $v$, and the other parameters are as the same as Example 10.1.

1) $g_{sm} = 0.001$

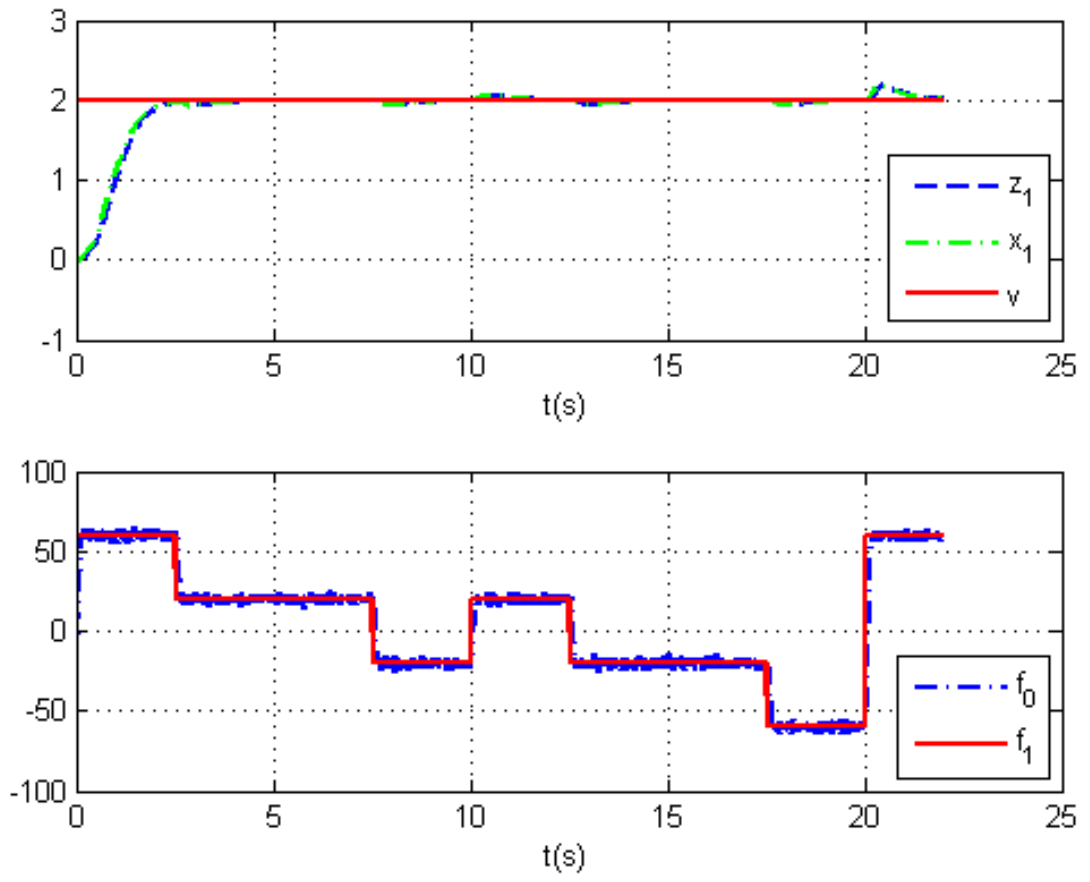

Fig.17. Step response ( $g_{sm} = 0.001$ )

2) $g_{sm} = 0.01$

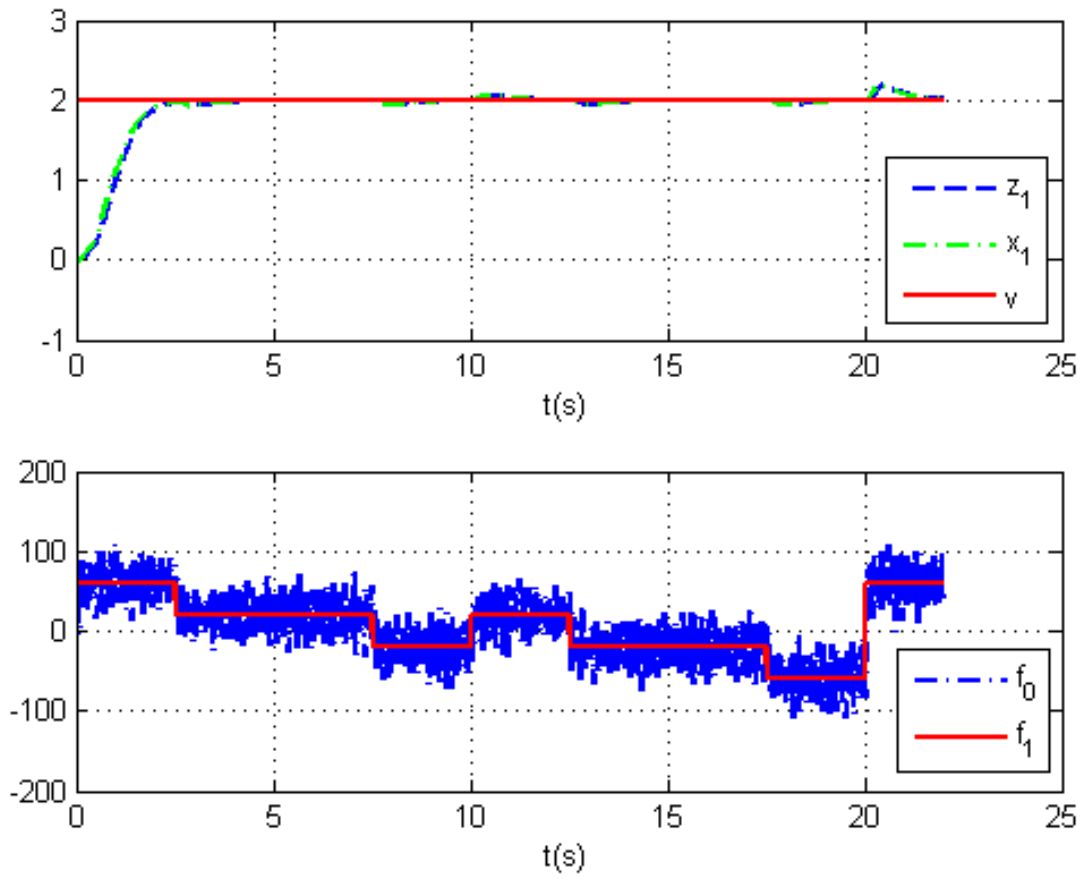

Fig.18. Step response ( $g_{sm} = 0.01$ )

3) $g_{sm} = 0.1$

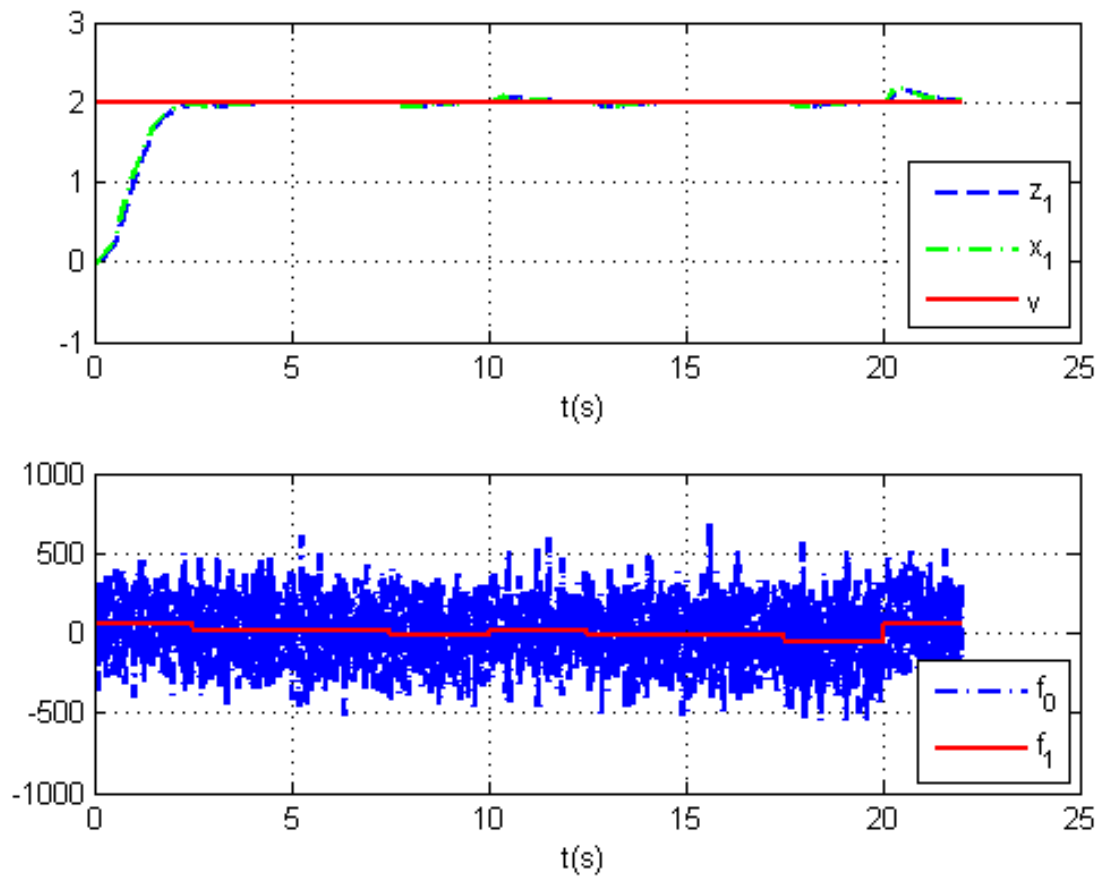

Fig.19. Step response ( $g_{sm} = 0.1$ )

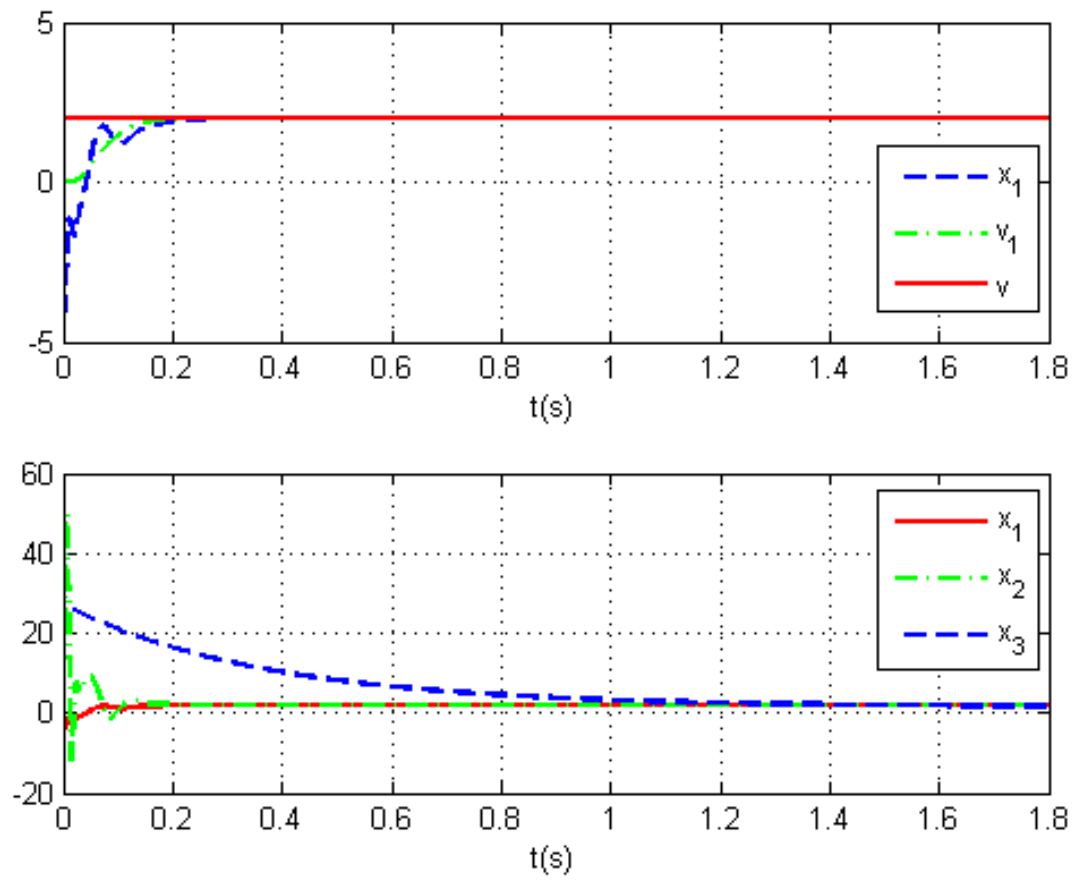

Fig.20. Step response ( $g_{sm} = 0.001$ )

When $f_1(t)$ is replaced by a staircase waveform, the estimation algorithm (241) can also well counteract the influence of disturbance, and has a strong ability to suppress noise.

**Example 10.3.** Lorenz chaotic control system (Han, 2008)

$$\begin{cases} Dx_1 = \sigma(x_2 - x_1) \\ Dx_2 = \rho x_1 - x_2 - x_1 x_3 + u \\ Dx_3 = -b_1 x_3 + x_1 x_2 \\ \chi = x_1 + v_m g_{sm1} wgn(1, length, -20) \end{cases} \tag{260}$$

Where, $\sigma = 10$, $\rho = 28$, $b_1 = \dfrac{8}{3}$, $x_1(0) = -4.47$, $x_2(0) = -0.505$, $x_3(0) = 28.02$.

(260) can be written as follows.

$$\begin{cases} D^2 x_1 = \sigma\left[\rho x_1 - x_2 - x_1 x_3 + \sigma(x_1 - x_2)\right] + \sigma u \\ Dx_3 = -b_1 x_3 + x_1 x_2 \\ \chi = x_1 + v_m g_{sm} wgn(1, length, -20) \end{cases} \tag{261}$$

So the Lorenz system (260) has a relative order of two.

The given signal is $v$, and the other parameters are as follows.

$v = v_m = 2$, $h = 0.0001s$; arranged transition process: $m_1 = 4$, $n_1 = 200$; ESO: $m_2 = 4$, $n_2 = 20$; controller: $m_3 = 2$, $n_3 = 200$

1) $g_{sm} = 0.001$

2) $g_{sm} = 0.01$

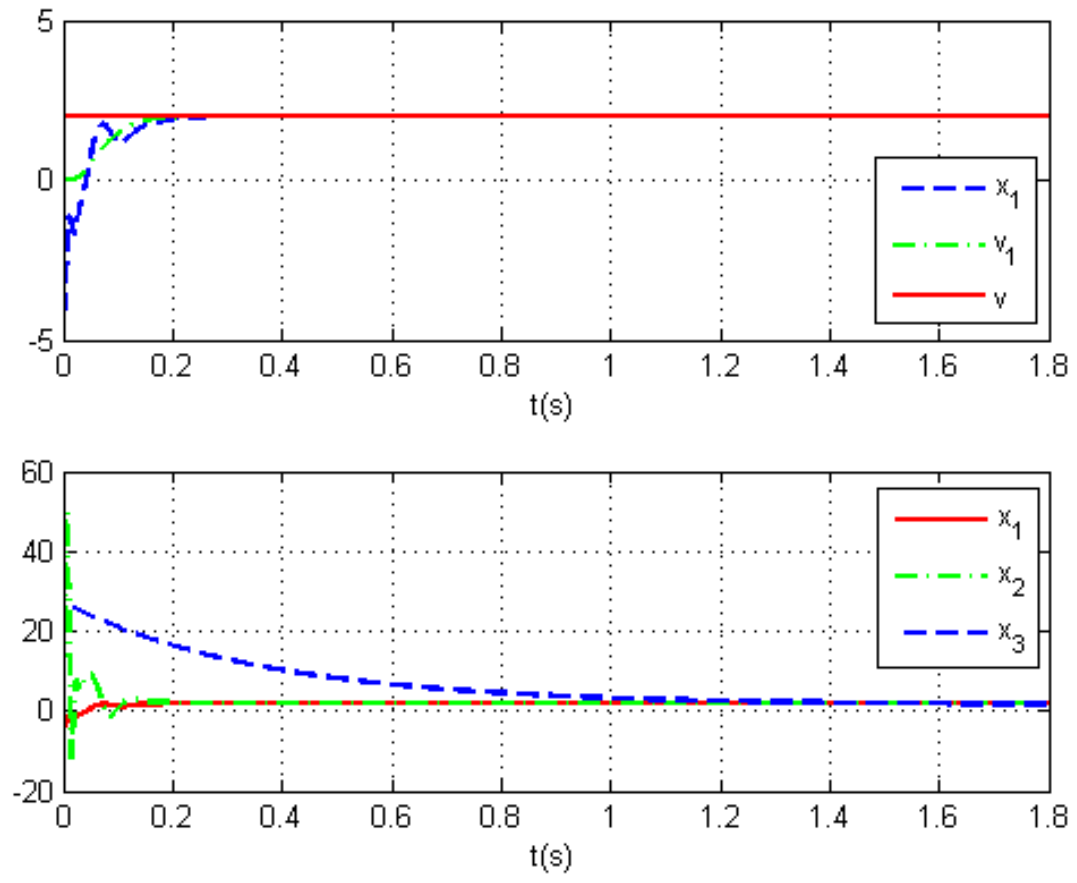

Fig.21. Step response ( $g_{sm} = 0.01$ )

3) $g_{sm} = 0.1$

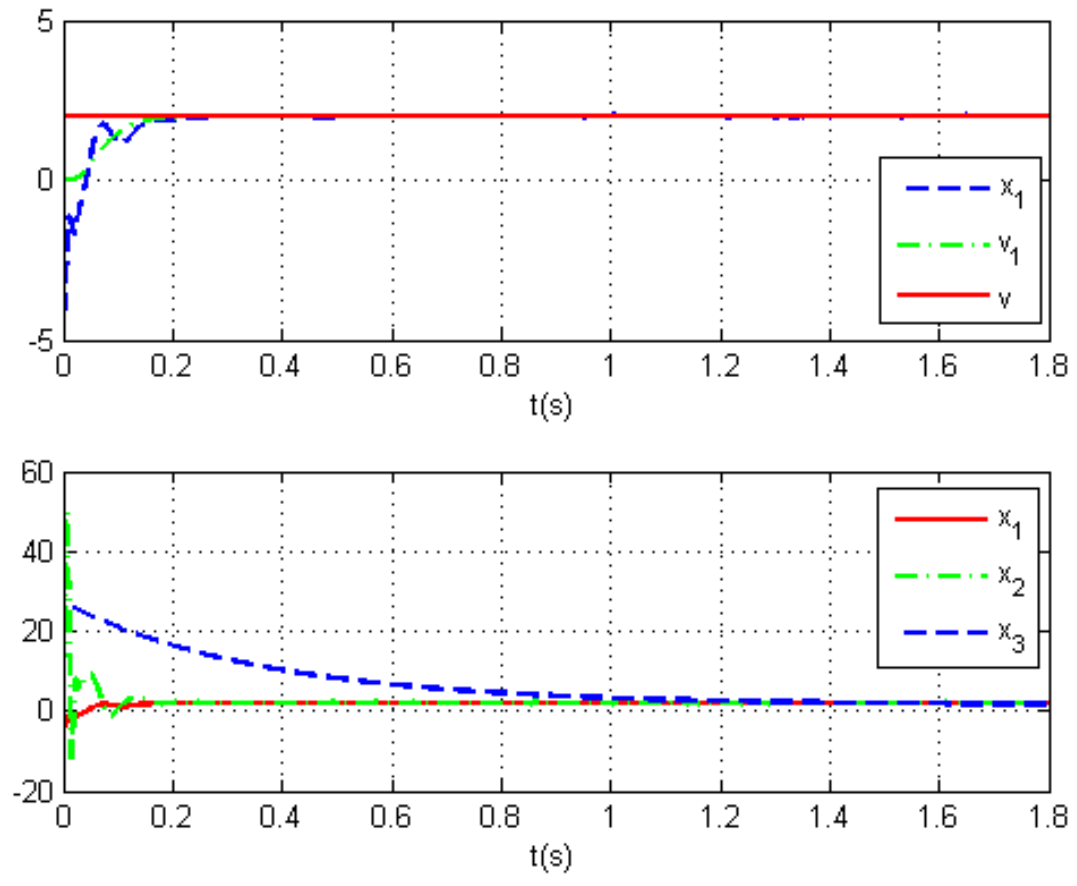

Fig.22. Step response ( $g_{sm} = 0.1$ )

$v_1$ is the arranged transition process. From Fig.20 to Fig.22, it can be seen that the magnitude of Gaussian white noise in output signal (or feedback signal) has almost no effect on the step response of Lorenz chaotic system, and the signals of $x_1$, $x_2$ and $x_3$ converge to the constants, respectively. This shows that the Lorenz chaotic system can be effectively controlled by the ADRCT which has a strong ability to suppress noise.

The above simulation and analysis show that the new ADRCT built in this paper has good immunity and robustness.

## 11. Conclusion

This paper first introduces the basic concepts of generating function in combinatorics, and then derives the 2nd-order synthesis function $fsun()$. By means of generating function, $fxiao()$, which is the m-order time optimal control synthesis function of discrete system with the origin as the end point, is derived in detail, and the m-order tracking-form synthesis function is given. In order to solve the contradiction between the rapid demand for time optimal control and the good noise immunity requirement for control system, a filter factor is introduced and only analyzed in the linear case. The relevant conclusions obtained are used to extract and predictively compensate the filtered feedback signal and its differential signals. The phase delay and amplitude attenuation of extracted signal and its differential signals are actually produced and verified by numerical simulation, and the numerical simulation shows that the predictive compensation method proposed in this paper can effectively compensate them.

The reconstruction of ADRCT is a successful application of m-order synthesis function in control science and engineering. Based on the m-order tracking-form synthesis function with filter factor, the new type structures of ESO and feedback controller that satisfy the separation principle are reconstructed to construct a new ADRCT, and the correctness of dynamic compensation linearization method used in this paper is proved. Finally, the numerical simulation shows that the ADRCT has good immunity and robustness. It can be expected that the new ADRCT has a wide range of applications in the field of control engineering.

## Acknowledgements

The author of this paper would like to extend his sincere thanks to the authors of (Qian & Song, 2011) and (Han, 2008) for their in-depth study of time optimal control system and the author of (Han, 2008) who proposed and developed the ADRC technology, and also would like to thank the authors of (Sun & Sun, 2010), who put forward the 2nd-order synthesis function, for their wonderful article that immediately led the author of this paper to understand the correct direction for handling the topic of this paper.